\documentclass{amsart}

\usepackage{amsfonts,amssymb, amsmath}
\usepackage{epsf}
\usepackage{graphicx}
  \newtheorem{theorem}{Theorem}[subsection]
\newtheorem{lemma}[theorem]{Lemma}
\newtheorem{corollary}[theorem]{Corollary}

\newtheorem{prop}{Proposition}[subsection]
\theoremstyle{definition}
\newtheorem{definition}[theorem]{Definition}

\theoremstyle{remark}
\newtheorem{remark}[theorem]{Remark}

\numberwithin{equation}{subsection}

\begin{document}

\title[Simple Approximately Subhomogeneous algebras]
{On the classification of simple approximately subhomogeneous C*-algebras not necessarily of real rank zero}
\author{C. Ivanescu}
\address{Department of Mathematics, University of Toronto,
100 St. George Street, Toronto, ON, M5S 3G3, Canada}
\email{cristian@math.toronto.edu }

\large
\begin{abstract} A classification is given of certain separable nuclear C*-algebras not necessarily of real rank zero, namely the class of simple C*-algebras which are inductive limits  of continuous-trace C*-algebras whose building blocks have their spectrum homeomorphic to the interval $[0,1]$ or to a finite disjoint union of closed intervals. In particular, a classification of those stably AI algebras which are inductive limits of hereditary sub-C*-algebras of interval algebras is obtained. Also, the range of the invariant is calculated. 
\end{abstract}
\maketitle

\section{Introduction}
 Q. Lin and N. C. Phillips have shown in \cite{lin} that an important class of simple crossed product C*-algebras are approximately subhomogeneous algebras, abbreviated ASH algebras. We recall that a C*-algebra is said to be subhomogeneous if it is isomorphic to a sub-C*-algebra of $pM_n(C_0(X))p$ for some natural number $n$, $p$ a projection in $M_n(C_0(X))$ and for some locally compact Hausdorff space $X$. An ASH algebra is an inductive limit of subhomogeneous algebras. \\
\indent It has become an important task to classify simple ASH algebras by their Elliott invariant. This article contains partial results in this direction.\\
\indent The first result on the classification of non real rank zero algebras was the classification by G. Elliott of unital simple approximately interval algebras, abbreviated AI algebras (see \cite{ell2}). This result was extended to the non-unital case independently by I. Stevens (\cite {ste}) and K. Thomsen (\cite {kto}). Also, an interesting partial result in the non-simple case was given by K. Stevens (\cite {kst}). It is worth mentioning that all these algebras are approximately homogeneous algebras, abbreviated AH algebras, and that the most general classification result for AH algebras was obtained by Elliott, Gong and Li in \cite {go2}.\\ 
\indent The first isomorphism result for ASH algebras was the proof by H. Su of the classification of C*-algebras of real rank zero which are inductive limits of matrix algebras over non-Hausdorff graphs; see \cite{su}. The classification of ASH algebras was also considered in \cite {go}, \cite {hli}, \cite {mig}, \cite{rak}, and \cite {kto2}. (The list of contributions is intended to be representative rather then complete for the classification of ASH algebras.)\\ 
\indent An important result on the classification of ASH algebras not of real rank zero, and in fact the first one, is due to I. Stevens (\cite{ste0}). The main result of the present article is a substantial extension of Stevens's classification to the class given by all simple C*-algebras which are inductive limits of continuous-trace C*-algebras with spectrum homeomorphic to the closed interval $[0,1]$ (or a finite disjoint union of closed intervals). In particular, the spectra of the building blocks considered here are the same as for those considered by Stevens. The building blocks themselves are more general.\\  
\indent We proceed by approximating the building blocks appearing in a given inductive limit decomposition with special building blocks. Here by special building blocks we mean continuous-trace C*-algebras with finite dimensional irreducible representations and such that the dimension, as a function on the interval, is a finite (lower semicontinuous) step function. We prove that in this case the algebra has a finite presentation, with stable relations. A result of T. Loring allows us to conclude that we can replace the original inductive limit decomposition with an inductive limit of special building block C*-algebras.\\  
\indent The advantage of having inductive limits of special subhomogeneous algebras is the following: for such special subhomogeneous algebras we have an Existence Theorem and a Uniqueness Theorem due to I. Stevens. Making use of this advantage we prove the classification theorem by applying the Elliott intertwining argument.\\
\indent  As a special case of the theorem we obtain the classification of simple inductive limits of arbitrary hereditary sub-C*-algebras of interval algebras.(Stevens considered very special inductive limits of hereditary sub-C*-algebras.) \\
\indent We also complete I. Stevens's description of the range of the invariant by including the case of an unbounded trace norm map.\\
\indent We end the article by comparing the size of the class of stably AI algebras and the class of AI algebras.

\emph{ \bf Acknowledgments.} I would like to thank George Elliott for his many helpful comments on earlier drafts of
this paper which led to a substantial increase in clarity.  This research was supported by both the Department of Mathematics, University of Toronto and by Ontario Graduate Scholarship.

\section{The problem and the results} 
 We consider the following class of C*-algebras: simple C*-algebras which are inductive limit of continuous-trace C*-algebras whose building blocks have their spectra homeomorphic to $[0,1]$ or finite disjoint union of closed intervals.\\
\indent We prove a complete isomorphism theorem for this class, namely:
\begin{theorem} Let ${\mathcal A}$ and ${\mathcal  B}$ be two simple C*-algebras which are inductive limit of continuous-trace C*-algebras whose spectrum is homeomorphic to $[0,1]$ or finite disjoint union of closed intervals and all irreducible representations are finite dimensional. We also assume that the connecting maps are injective. Assume that \\
 1. there is an isomorphism $ \psi _0 : D({\mathcal A}) \rightarrow D({\mathcal B}),$\\
 2. there is an isomorphism $ \psi _T : (\mathrm{Aff}T^+{\mathcal A},\mathrm{Aff}'{\mathcal A}) \rightarrow (\mathrm{Aff}T^+{\mathcal B},\mathrm{Aff}'{\mathcal B})$\\
 and the two isomorphisms are compatible: 
 $$ \widehat{ \psi_0([p])} = \psi_T( \widehat{[p]}),\;\forall [p] \in D({\mathcal A}).$$\\
 Then there is an isomorphism of the algebras ${\mathcal A}$ and ${\mathcal B}$ that induces the given isomorphism at the level of the invariant.
\end{theorem}
\begin{corollary} Let ${\mathcal A}$ and ${\mathcal B}$ be two simple C*-algebras which are inductive limit of hereditary sub-C*-algebras of interval algebras. We also assume that the connecting maps are injective. Assume that \\
 1. there is an isomorphism $ \psi _0 : D({\mathcal A}) \rightarrow D({\mathcal B}),$\\
 2. there is an isomorphism $ \psi _T : (\mathrm{Aff}T^+{\mathcal A},\mathrm{Aff}'{\mathcal A}) \rightarrow (\mathrm{Aff}T^+{\mathcal B},\mathrm{Aff}'{\mathcal B})$\\
 and the two isomorphisms are compatible: 
 $$ \widehat{ \psi_0([p])} = \psi_T( \widehat{[p]}),\;\forall [p] \in D({\mathcal A}).$$\\
 Then there is an isomorphism of the algebras ${\mathcal A}$ and ${\mathcal B}$ that induces the given isomorphism at the level of the invariant.

\end{corollary}
 We also describe the range of the invariant. More precisely, we prove the following theorem:
\begin{theorem} Suppose that $G$ is a simple countable dimension group, $V$ is the  cone associated to a metrizable Choquet simplex $S$ , $\lambda : S \rightarrow Hom^+(G,R)$ is a continuous affine map which sends extreme rays to extreme rays, and $f:S\rightarrow (0,\infty]$ any affine lower semicontinuous map. Then $[G,(V,S), \lambda,f]$ is the Elliott invariant of some simple non-unital inductive limit of continuous-trace C*-algebras.
\end{theorem} 

\section{Partial case considered by I. Stevens}
 I. Stevens \cite{ste0} has solved the classification problem in the particular case of algebras that are inductive limits of special building blocks: 
\begin{theorem}({\bf Stevens's Isomorphism Theorem}) Let ${\mathcal A}$ and ${\mathcal B}$ be simple approximately step-hereditary C*-algebras with injective connecting maps. Assume that \\
 1. there is an isomorphism $ \psi _0 : D({\mathcal A}) \rightarrow D({\mathcal B}),$\\
 2. there is an isomorphism $ \psi _T : (AffT^+{\mathcal A},Aff'{\mathcal A}) \rightarrow (AffT^+{\mathcal B},Aff'{\mathcal B})$\\
 and the two isomorphisms are compatible: 
 $$ \widehat{ \psi_0([p])} = \psi_T( \widehat{[p]}),\;\forall [p] \in D({\mathcal A}).$$\\
 Then there is an isomorphism of the algebras ${\mathcal A}$ and ${\mathcal B}$ that induces the isomorphism at the level of the invariant.
\end{theorem}
 I. Stevens \cite{ste0} also describes the range of the invariant, assuming the trace norm map is bounded:
\begin{theorem}({\bf Stevens's Range of the Invariant Theorem}) Let  $(G,u)$ be  a simple dimension group, $u$ a positive element in $G$, $(\Gamma,S)$ a cone with a metrizable Choquet simplex $S$ as its base, $\mu :S \rightarrow (0,1]$, a lower semicontinuous affine map and  $\lambda : \Gamma \rightarrow s(G)$, where $s(G)$ are the functionals on $G$, and $\lambda $ is a continuous, surjective affine map  sending  extreme rays in extreme rays. Then $[G,(V,S), \lambda,f]$ is the Elliott invariant of some hereditary subalgebra of a simple unital AI-algebra.
\end{theorem} 
\indent Our goal in this paper is to prove Theorem 2.0.1 and Theorem 2.0.3 which generalize Stevens's theorems 3.0.4 and 3.0.5.

\section{ Building blocks and their inductive limits}
\indent We start by  describing the class of algebras in which we are interested. Suppose ${\mathcal A}$ is a C*-algebra with Hausdorff spectrum $T$. It is known that each primitive quotient $A(t)$ has up to equivalence a unique irreducible representation $\pi _{t}$. Hence, whenever $p \in {\mathcal A}$ is an element of ${\mathcal A}$ such that $p(t)$ is a projection, the rank of $p(t)$ is well-defined as the dimension of the range of $\pi_{t}(p(t))$. So it make sense to say that $p(t)$ is a rank one projection.
\begin{definition} A continuous-trace C*-algebra is a C*-algebra ${\mathcal A}$ with Hausdorff spectrum $T$ such that, for each $t \in T$, there are a neighborhood $V$ of $t$ and $a \in {\mathcal A}$ such that $a(s)$ is a rank one projection for all $s \in V$.
\end{definition}
\begin{remark} The existence of a local rank one projection is sometimes summed up by saying that ${\mathcal A}$ $satisfies$ $Fell's$ $condition$.
\end{remark}
\begin{remark} If ${\mathcal A}$ is a continuous-trace C*-algebra, then its spectrum $T$ is a locally compact Hausdorff space. By \cite {dix0} we get a continuous field $\mathcal A$ over $T$ of elementary C*-algebras with Fell's condition, and $A \cong \Gamma _{0}( \mathcal{A})$, the C*-algebra of cross-sections, which maps $C(T)$ into the center $Z(M({\mathcal A}))$ of the multiplier C*-algebra $M({\mathcal A})$.
\end{remark}
\indent We recall:
\begin{definition} A C*-algebra ${\mathcal A}$ is stable if ${\mathcal A} \otimes \mathbb{K}\cong {\mathcal A}$, with $\mathbb{K}=\mathbb{K}(H)$ the algebra of compact operators on a separable infinite dimensional Hilbert space $H$.
\end{definition}

\begin{remark} A stable C*-algebra is necessarily a non-unital C*-algebra.
\end{remark}

\begin{definition} Two C*-algebras ${\mathcal A}$ and ${\mathcal B}$ are stably isomorphic, denoted  ${\mathcal A} \sim {\mathcal B}$, if ${\mathcal A} \otimes \mathbb{K }\cong {\mathcal B} \otimes \mathbb{K}$. 
\end{definition}

\begin{remark}The notion of being stably isomorphic for C*-algebras (or Morita-Rieffel equivalence for the separable case, cf. \cite {rif}) is an equivalence relation which is weaker than the isomorphism relation. It is easy to see that $M_{2}(\mathbb{C})$ is stably isomorphic to $M_{3}(\mathbb{C})$ but $M_{2}(\mathbb{C})$ is not isomorphic to $ M_{3}(\mathbb{C})$.
\end{remark}

\begin{definition} A C*-algebra ${\mathcal A}$ is a stably (isomorphic) AI algebra if ${\mathcal A} \otimes \mathbb{K}$ is an AI-algebra.
\end{definition}
\indent In a similar way we define:
\begin{definition} A C*-algebra is a stably (isomorphic) AF algebra (or AT algebra, or AH algebra) if ${\mathcal A}\otimes \mathbb{K}$ is an AF-algebra (or AT-algebra, or AH-algebra.)
\end{definition}
\begin{remark} Any AI algebra is also a stably AI algebra.
\end{remark}
\begin{remark} Any full hereditary subalgebra of an AI-algebra is a stably AI algebra. The class of algebras that Theorem 3.0.4 classifies, namely approximately step-hereditary subalgebras of simple AI-algebras, are stably AI algebras.
\end{remark}
\begin{remark} It can be shown (see \cite{dix0}) that in general a continuous-trace C*-algebra $A$ with spectrum $T$ is locally Morita-Rieffel equivalent to $C_{0}(T)$. The only obstruction for $A$ to be Morita-Rieffel equivalent to $C_{0}(T)$ is the Dixmier-Duoady class: $\delta(A)$ in $H^{3}(T,\mathbb{Z})$. In the case we are considering in this article, $T$ is $[0,1]$ or finite disjoint union of closed intervals for which the Dixmier-Duoady class vanishes. Therefore our building blocks are stably AI algebras, in fact stably interval algebras.
\end{remark}
\indent Next we notice that there are algebras which are stably AI algebras but not AI-algebras as opposed to stably AF algebras which are necessarily AF-algebra:
\begin{prop} Any stably AF algebra is necessarily an AF-algebra.
\end{prop}
\begin{proof}
 Let $A$ be a stably AF-algebra. $A \otimes \mathbb{K}$ is AF-algebra. If we denote by $e_{11}$ the rank one projection and by $a$ the unit of the multiplier algebra of $A$ then $A$ is isomorphic to the cut down $a \otimes e_{11}(A \otimes \mathbb{K})a \otimes e_{11}$ which is a hereditary subalgebra of the AF-algebra $A \otimes \mathbb{K}$. By Elliott's result \cite{ell1} we know that hereditary subalgebras of AF-algebras are AF-algebras and so $A$ is an AF-algebra.
\end{proof}
\begin{prop}({\bf I. Stevens}) There are algebras which are stably AI algebras but not AI-algebras.
\end{prop}  
\begin{proof}
 This follows from the existence of a hereditary subalgebra of a simple AI-algebra which is not AI-algebra. To prove this suppose the contrary. Then any hereditary subalgebra of an AI-algebra has an approximate unit consisting of projections and thus any AI-algebra has real rank zero. This is a contradiction because the AI-algebra $C[0,1]$ is not of real rank zero. The first example of a simple AI algebra not of real rank zero appeared in \cite{blab}. Uncountable many examples of simple AI-algebra not of real rank zero were constructed by K.Thomsen in \cite{kto}.
\end{proof}
 Moreover we have the following result:
\begin{prop} Any (simple) stably AI algebra can be realized as a hereditary subalgebra of a (simple) AI-algebra.
\end{prop}
\begin{proof}
 Assume $A$ is a (simple) stably AI-algebra. Hence $A \otimes \mathbb{K}$ is a (simple) AI-algebra. $A$ is isomorphic to the cut down: $$A \cong a \otimes e_{11}(A \otimes \mathbb{K})a \otimes e_{11}$$
 which is a hereditary subalgebra of the (simple) AI algebra $A \otimes \mathbb{K}$, where $a$ is the unit of the multiplier algebra of $A$, $e_{11}$ a rank $1$ projection in $\mathbb{K}$.
\end{proof}
\begin{remark} Any continuous-trace C*-algebra whose spectrum is homeomorphic to the closed interval $[0,1]$ and has all irreducible representations finite dimensional is isomorphic to a hereditary sub-C*-algebra of an interval algebra. 
\end{remark}
\indent A natural question to ask is whether all simple ASH algebras are stably AH algebras. This question has a negative answer:
\begin{prop} There are simple ASH algebras that are not stably AH algebras.
\end{prop}
\begin{proof} For simple inductive limit of splitting interval algebras is shown in \cite{su} that $K_0$ fails to have the Riesz decomposition property. Therefore these ASH algebras can not be stably AH algebras.
\end{proof}  

 It is important to notice that the class of (simple) stably AI algebras is closed under taking inductive limits: 
\begin{prop} Inductive limits of (simple) stably AI algebras are still (simple) stably AI algebra.
\end{prop}
\begin{proof}
 Consider an inductive sequence of (simple) stably AI algebras $(A_i,\psi_i)$. Then we know that $A_i \otimes \mathbb{K}$ is a (simple) AI algebra. Moreover the inductive limit of the direct system $(A_i \otimes \mathbb{K}, \psi \otimes id) $ 
 $$\lim\limits_{\rightarrow}(A_i \otimes \mathbb{K}) \cong (\lim\limits_{\rightarrow}A_i) \otimes \mathbb{K} $$
 is a simple AI-algebra by using \cite{ell3}, Theorem 4.3 the implication (iii)$\rightarrow$ (iv). Therefore $\lim\limits_{\rightarrow}A_i$ is (simple) stably AI algebras.
\end{proof}

 We end this section with an example of a continuous-trace C*-algebra $\mathcal{A}$ which is stably AI algebra but not an AI algebra. In addition $\mathcal{A}$ is a C*-algebra not of real rank zero.
$$\mathcal{A}= \left( \begin{array}{cc}
 C_{0}[0,1) &  C_{0}[0,1) \\
 C_{0}[0,1)  &  C[0,1]     \\
\end{array} \right) \subset M_{2}(C[0,1])$$
 It is not hard to check that $\mathcal{A}$ is Morita-Rieffel equivalent to $C[0,1]$. Since the real rank of $C[0,1]$ is not zero it follows that the real rank of $ \mathcal{A}$ is not zero. Moreover, $\mathcal{A}$ does not have an approximate unit consisting of projections. Therefore it can not be an AI algebra. In particular $\mathcal{A}$ is an example of a subhomogeneous C*-algebra which is not homogeneous algebra.

\section{Stably isomorphic algebras}

 Using the complete classification result of I. Stevens \cite {ste} for simple AI-algebras we show that the classification of simple stably AI algebras reduces to the classification of hereditary subalgebras of simple AI algebras which are stably isomorphic. 

\begin{prop} Let $H_1$ and $H_2$ be two hereditary subalgebras of simple AI-algebras  $A_1$ and respectively $A_2$. Then if $A_1$ and $A_2$ are not stably isomorphic then $H_1$ is not isomorphic to $H_2$. On the contrary: when $A_1$ is stably isomorphic to $A_2$ we obtain that $H_1$ is stably isomorphic to $H_2$.
\end{prop}
\begin{proof}
 Because of simplicity of $A_1$ and $A_2$  we get that the hereditary subalgebras $H_1$ and $H_2$ are full. Using Brown's result \cite {bro} we obtain that $H_1$ is stably isomorphic to $A_1$ and  $H_2$ is stably isomorphic to $A_2$. As a consequence, if $A_1$ is not stably isomorphic to $A_2$ then $H_1$ is not isomorphic to $H_2$. In the other case, when  $A_1$ is stably isomorphic to $A_2$ we get that $H_1$ is stably isomorphic to $H_2$.
\end{proof} 

\begin{remark}
 However, to prove the isomorphism theorem we will not restrict our attention to the class of simple inductive limits of continuous-trace C*-algebras which are Morita-Rieffel equivalent.
\end{remark}

\section{Spectrum of the building blocks}

\indent We start by noticing that the spectrum of a full hereditary sub-C*-algebra  of an interval algebra is a compact subset of the spectrum of the interval algebra. Here by full we mean that there are no non-trivial ideals of the interval algebra which contain the hereditary sub-C*-algebra. Given a sub-C*-algebra of a simple unital AI algebra which is the inductive limit of hereditary sub-C*-algebras of interval algebras, in \cite{ste0} was shown that almost all but finitely many hereditary subalgebras are full, in particular their spectrum is a closed subset of the spectrum of the interval algebra. As a more general result we prove:

\begin{theorem} Let $A=\lim\limits_{\rightarrow}A_i$ be a simple C*-algebra, where the maps are assumed to be injective and $A_i$'s are assumed to be hereditary sub-C*-algebras of interval algebras. Then $A$ can be realized as the inductive limit of hereditary sub-C*-algebras of interval algebras such that all but finitely many building blocks have compact spectrum. More precisely the sequence can be chosen with the property that each building block contains a full projection. 
\end{theorem}

\begin{proof}
 Since $A$ is a simple C*-algebra, by [\cite{blc}, Corollary 5.2] we obtain that $A\otimes O_2$ contains a non-trivial projection $p$. Moreover we have $$A\otimes O_2 = \lim\limits_{\rightarrow}(A_i\otimes O_2),$$
and the spectrum of the building blocks remains the same up to isomorphism, $\widehat{A_i \otimes O_2} \cong \widehat{A_i}$.\\
\indent Next we pull back the projection $p$ at finite stages. Hence we can assume that $p \in A_i\otimes O_2$ for $i\geq j_0$.\\
\indent Consider the ideal generated by $p$ inside $A_i \otimes O_2$ which we denote by $I_i$. Then $I_i$ is a non-zero ideal and it corresponds to a compact subset of the spectrum of $A_i \otimes O_2$.\\
\indent Let $I=\lim\limits_{\rightarrow}I_i$. Then $I$ is a non-zero ideal of $A \otimes O_2$ which is a simple C*-algebra. Therefore $$I=A \otimes O_2.$$
\indent Because $O_2$ is a simple C*-algebra we have $I_i \cong \hat{I_i} \otimes O_2$ with $\hat{I_i}$ ideal of $A_i$. Therefore $$A= \lim\limits_{\rightarrow}\hat{I_i}$$ where each $\hat{I_i}$ is a closed two-sided ideal of the given building block $A_i$. Moreover by construction the new building blocks $(\widehat{I_i})_i$ contain a full projection and they have compact spectrum.
\end{proof}
\begin{remark} A similar proof can be obtained by noticing that $A \otimes \mathbb{K}$ is an AI algebra which contains a non-trivial projection.
\end{remark}
\begin{remark}
 As a consequence of the above theorem, in the isomorphism theorem 2.0.1 we can allow inductive limits of continuous-trace C*-algebras whose spectrum is any locally compact subset of the closed interval $[0,1]$.
\end{remark}
\begin{remark}  
We also note that for the isomorphism theorem 2.0.1 there is no loss of generality in assuming that the building blocks are full hereditary sub-C*-algebra of interval algebra.
\end{remark}
 
\section{The Dimension Function} 
  A very important data that we consider is a map that assigns to each class of irreducible representations, the dimension of a representation from that class. 
More precisely pick a Hilbert space $H$ such that every irreducible *-representation can be realized in a closed subspace of $H$.

\begin{definition} Let $A$ be a C*-algebra and with $\hat{A}$ denote the spectrum of $A$. Then the $dimension$ $function$ is a map from $\hat{A}$ to $ \mathbb{R} \cup \infty $:
            $$ \pi \rightarrow dim(H_{\pi}),  $$ where by $dim(H_{\pi})$ we mean the dimension of the irreducible representation $\pi$.
\end{definition}
\begin{remark} 
 A closely related definition appears in (\cite{ped0}, pp.197).
\end{remark}
\begin{prop} The dimension function is a lower semicontinuous function.
\end{prop}
\begin{proof} See for instance \cite{ped0}, pp.73. 
\end{proof}

\begin{remark}
 The dimension function of a commutative algebra takes values in the set $\{ 0, 1 \}$.
\end{remark}
\begin{remark} Let $d:X \rightarrow \mathbb{Z}$ be a lower semicontinuous integer-valued function, where $X$ is a locally compact space. Then in [\cite{vas}, Theorem 3 and Theorem 4] it is shown by N. B. Vasil'ev how to describe a subhomogeneous C*-algebra whose dimension function is precisely the function $d$.
\end{remark}

\begin{remark} Let $A$ be a subhomogeneous C*-algebra whose dimension function is $d:\hat{A} \rightarrow \mathbb{Z}$. There exists a canonical projection-valued function that gives rise to the dimension function $d$ obtained as follows. Let $H$ be a Hilbert space such that every irreducible *-representation can be realized in a closed subspace of $H$. Next define $$P^{A}: \hat{A} \rightarrow {\bf B}(H),$$ 
 $$P^{A}(\pi_{0})=I_{\pi_{0}}$$ where $\pi_0 \in \hat{A}$ and by $I_{\pi_{0}}$ we denote the unit of ${\bf B}(H_{\pi_{0}})$ which is the canonical projection onto the Hilbert subspace $H_{\pi_0}$.\\
\indent The projection-valued function $P^{A}$ is a lower semicontinuous function in the following sense: $$[0,1]\ni \hat{\pi} \rightarrow <P^{A}(\hat{\pi})v|v>$$ is a real-valued lower semicontinuous function, where $v$ is a vector in the Hilbert space $H$.  
\end{remark}

\begin{remark} Let $P:[0,1] \rightarrow {\bf B}(H)$ be a lower semicontinuous projection-valued function such that the maximal rank of the projections appearing in the range of $P$ is $n$. Inside the full matrix algebra $M_n(\mathbb{C}) \otimes C[0,1]$ take all continuous matrix-valued functions which are smaller than $P$ and consider the hereditary sub-C*-algebra of $M_n(\mathbb{C}) \otimes C[0,1]$ that they generate. In this manner we have constructed a hereditary sub-C*-algebra which has the dimension function precisely the given $P$ composed with rank function.
\end{remark}

\subsection { Special subalgebras of interval algebras} Following a terminology introduced in \cite {ste0} we define:
\begin{definition} A continuous-trace C*-algebra whose spectrum is $[0,1]$ (or a hereditary sub-$C^*$-algebra of an interval algebra $C[0,1] \otimes F$) is a  special continuous-trace C*-algebra (respectively $step$ $hereditary$ sub-C*-algebra) if its dimension function is a $finite$ $step$ function: there is a partition of $[0,1]$ into a finite union of intervals such that the dimension function is constant on each such subinterval.
\end{definition}
\begin{definition} A hereditary sub-C*-algebra of an interval algebra is called step hereditary (see \cite{ste0}) if the unit of the bidual $p$ is a step projection: there is a partition of $[0,1]$ into a finite union of intervals such that $p$ is the same projection on each such subinterval.
\end{definition}

\begin{remark} A priori our definition for a step hereditary sub-C*-algebra is more general than I. Stevens's definition. As we will show in this paper, any step hereditary sub-C*-algebra in our sense is isomorphic to a step hereditary sub-C*-algebra in I. Stevens's sense.
\end{remark}

\section{ An isomorphism result for the building blocks}

\indent In this section we prove that for a certain class of continuous-trace C*-algebras, the dimension function is a complete invariant.\\
\indent We start by showing several lemmas that we use in proving the main result of this section:

 \begin{lemma} Let $f:[0,1] \rightarrow F$ be a lower semicontinuous projection-valued function, where $F$ is a simple finite dimensional C*-algebra. Assume that all projections appearing in the range of $f$ have the same rank, namely equal to $k$. Then $f$ is a continuous projection-valued function.
\end{lemma}

\begin{proof} Consider an arbitrary point $t_0$ in $[0,1]$ and a sequence $(t_n)_n$ in the interval $[0,1]$ that converge to $t_0$.\\
\indent Since $(f(t_n))_n$ is a sequence inside of the compact set of projections of a fixed rank $k$ we know that there exists a convergent subsequence. Without loss of generality we assume
 $$(f(t_n))_{n\in \mathbb{N}} \rightarrow e,$$
 where $e$ is a projection of rank $k$. \\
\indent But $f$ is a lower semicontinuous function and hence
 $$f(t_0) \geq e.$$
\indent Notice that $f(t_0)$ is a projection whose rank is the same as the rank of $e$. It is known that any two comparable projections with the same rank in a simple finite dimensional C*-algebra must be equal. Therefore:
$$\mathrm{lim}f(t_n)= f(t_0)=e.$$
\indent Hence the function $f$ is continuous.
\end{proof}

\begin{lemma} Let $f:[0,1] \rightarrow \mathbb{K}$ a lower semicontinuous projection valued function, where $\mathbb{K}$ is the algebra of compact operators. Assume that all projections appearing in the range of $f$ have the same finite rank $k$. Then there is a continuous function $V:[0,1] \rightarrow \mathbb{K}$ with values partial isometries and a projection $p \in \mathbb{K}$ such that 
 $$V(t)V(t)^{*} = p ,\;\;\;\forall t \in [0,1],$$
 $$V(t)^{*}V(t) = f(t),\;\;\;\forall t \in [0,1].$$
\end{lemma}

\begin{proof}
 Consider a point $t_0$ in $[0,1]$ and $p$ a projection of rank $k$. Following ideas of Glimm from \cite{gli} for each $t_0 \in [0,1]$ we can find $V(t)\in C[0,1]\otimes \mathbb{K}$ such that $V(t_{0})V^{*}(t_{0})=p$, $V^{*}(t_{0})V(t_{0})=f(t_{0})$ and $||V(t)p(t)V(t)^{*}-f(t)||< \frac{1}{2}$ for all $t$ in a neighborhood of $t_0$.\\
\indent Making use of a continuous function $g:[0,1]\rightarrow [0,2]$ defined to be zero in a neighborhood of zero and $\frac{1}{t}$ for $t\geq \frac{1}{2}$ we can define a continuous family of partial isometries $$W(t)=p[g(p(t)V(t)f(t)V(t)^{*}p]^{1/2}V(t)f(t)V(t)^{*}$$ such that for all $t$ in a neighborhood of $t_0$ we have: $$W(t)W^{*}(t)=p(t),$$ $$W^{*}(t)W(t)=f(t).$$
\indent By compactness of $K$ we find a finite open cover such that on each open set there is a continuous family of partial isometries that realize the equivalence of projections $f(t)$ and $p$ and with any open set overlapping with at least another open set. Next we proceed as in [\cite{kto1}, Lemma 2 or \cite{ell2}] to glue together the families of partial isometries on the overlapping sets. Consequently we obtain the desired family of partial isometries on the connected compact set $K$. By letting $K=[0,1]$ we complete the proof of the lemma.
\end{proof}

\begin{corollary} Let $f:[0,1] \rightarrow M_n(\mathbb{C})$ a lower semicontinuous projection valued function. Assume that all projections appearing in the range of $f$ have the same finite rank $k$. Then there is a continuous function $U:[0,1] \rightarrow M_n(\mathbb{C})$ with values unitaries and a projection $p \in M_n(\mathbb{C})$ such that: 
 $$U(t)U(t)^{*} = p ,\;\;\;\forall t \in [0,1],$$
 $$U(t)^{*}U(t) = f(t),\;\;\;\forall t \in [0,1].$$
\end{corollary}

\begin{proof}

 Consider a point $t_0$ in $[0,1]$. We proceed as in the previous lemma and find an open set $V_0$ containing $t_0$ and a continuous family of unitaries $U(t)_{t \in [0,1]}$ such that:
  $$U_{0}(t)f(t)U_{0}(t)^{*}=p_0,\;\;\;\forall t \in V.$$
\indent It is possible to find a continuous family of unitaries because in $M_n(\mathbb{C})$ if $p$ and $q$ are two Murray-von Neumann equivalent projections then $1-p$ and $1-q$ are another two Murray-von Neumann equivalent projections.\\ 
\indent By compactness of $[0,1]$ we find a finite open cover with any open set overlapping with at least another open set and having a similar property as $U_0$ has above:
$$ K \subset \bigcup_{i=1}^{N}V_i $$
\indent Given $i$ there is $j \neq i$ such that: $V_i \cap V_j \neq \oslash $ and 
 $$U_{i}(t)f(t)U_{i}(t)^{*}=p_i,\;\;\;\forall t \in V_i $$
\indent To get the desired family of unitaries we ``glue'' together the unitaries on  the overlapping set as follows:\\
\indent Let $t_1 \in V_i \cap V_j $. If $V_i$ is on the left side of $V_j$ then we modify $U_j$:
 $$ U^{ij}(t)  = \left\{ \begin{array}{cc}
                    U_i(t), & t \leq t_1\\
                    U_{i}(t_{1})U_{j}(t_1)^{*}U_{j}(t), & t > t_1
                   \end{array}
            \right. $$
\indent Notice that $p$ is equal to the projection corresponding to the most left $V_i$. For instance $$U^{ij}(t)f(t)(U^{ij}(t))^{*}=p_i,\;\;\mathrm{if}\;i\leq j\;\mathrm{and}\;t\in V_i\cup V_j,\;V_i \cap V_j \neq \varnothing.$$
\end{proof}

\begin{lemma} Let $f:[0,1] \rightarrow \mathbb{K}$ be a lower semicontinuous projection-valued function, where $\mathbb{K}$ is the algebra of compact operators. Assume that $f$ has constant rank $k$ on $[a,b)\subset [0,1]$. Then there is a continuous function $V:[a,b)\rightarrow \mathbb{K}$ with values partial isometries and a projection $p\in \mathbb{K}$ of rank $k$ such that:
 $$V(t)(V(t))^{*}=p\;\;for\;all\;t\in [a,b),$$
 $$(V(t))^{*}V(t)=f(t)\;\;for\;all\;t\in [a,b).$$
\end{lemma}

\begin{proof} The idea is to approximate the half-open interval with an increasing union of closed (hence compact) subsets. For example: 
 $$[a,b)= \bigcup_{n=1}^{\infty}[a,b-1/n].$$
\indent By Lemma 8.0.4 on each closed subset we have a continuous family of partial isometries. To extend from a smaller closed interval to a larger one we proceed as in [\cite {gli} or \cite{kto1}, Lemma 2 or \cite {ell2}].
\end{proof}

 \begin{corollary} Let $f:[0,1] \rightarrow M_n(\mathbb{C})$ be a lower semicontinuous projection valued function. Assume that $f$ has constant rank $k$ on $[a,b)\subset [0,1]$. Then there is a continuous function $U:[a,b)\rightarrow M_n(\mathbb{C})$ with unitaries values and a projection $p\in M_n(\mathbb{C})$ of rank $k$ such that:
 $$U(t)(U(t))^{*}=p\;\;for\;all\;t\in [a,b)$$
 $$(U(t))^{*}U(t)=f(t)\;\;for\;all\;t\in [a,b).$$
\end{corollary}

\begin{proof} As in the previous lemma we use 
  $$[a,b)= \bigcup_{n=1}^{\infty}[a,b-1/n]$$
By the Corollary 8.0.5 we have families of unitaries on each set of the form $[a,b-1/n]$. It remains to check that we can extend continuously this family of unitaries from the smaller compact subset to a larger one. For instance consider two compact sets $F_1$ and $F_2$ with $F_{1} \subset F_{2}$ and their corresponding unitaries $U_{1}(t),\;t \in F_1$ and respectively $U_{2}(t),\;t \in F_{2}$. \\
\indent Observe that we can order these intervals. Consider $t_0$ in $F_1$ not far from the boundary with $F_2$.
 Define $$ U^{1,2}(t)  = \left\{ \begin{array}{cc}
                    U_1(t) & t \in F_1\;\mathrm{and}\; t \geq  t_0\\\
                    U_{1}(t_{0})U_{2}(t_0)^{*}U_{2}(t) & t < t_0
                   \end{array}
            \right. $$
 and hence we get a continuous family of unitaries that have the desired property of the corollary on the set $[a,b)$.
\end{proof}

\begin{lemma} Let $f:[0,1] \rightarrow \mathbb{K}$ be a lower semicontinuous projection valued map, where $\mathbb{K}$ is the algebra of compact operators. Assume that $f$ has constant rank equal to $k$ on $(a,b)\subset [0,1]$. Then there is a continuous map $V:(a,b)\rightarrow \mathbb{K}$ with values partial isometries and a projection $p\in F$ of rank $k$ such that:
 $$V(t)(V(t))^{*}=p\;\;for\;all\;t\in (a,b)$$ 
 $$(V(t))^{*}V(t)=f(t)\;\;for\;all\;t\in (a,b).$$
\end{lemma}
\begin{proof} We have  $$(a,b)= \bigcup_{n=1}^{\infty}[a+1/n,b-1/n]$$ 
 and proceed as in the half open interval case, see Lemma 8.0.6, this time extending the family of partial isometries in both directions.
\end{proof}

 \begin{corollary} Let $f:[0,1] \rightarrow M_n(\mathbb{C})$ be a lower semicontinuous projection valued map. Assume that $f$ has constant rank $k$ on $(a,b)\subset [0,1]$. Then there is a continuous map $U:(a,b)\rightarrow M_n(\mathbb{C})$ with values unitaries and a projection $p\in M_n(\mathbb{C})$ of rank $k$ such that:
 $$U(t)(U(t))^{*}=p\;\;for\;all\;t\in (a,b),$$ 
 $$(U(t))^{*}U(t)=f(t)\;\;for\;all\;t\in (a,b).$$
\end{corollary}
\subsection {A decomposition result for the projection-valued functions}
 In Remark 7.0.27 we associated a lower semicontinuous projection-valued function to the dimension function of any subhomogeneous C*-algebra.\\
\indent Next we prove that there is a decomposition of such lower semicontinuous projection-valued functions as a sum of simpler lower semicontinuous projection-valued functions. As a consequence of this decomposition we prove that certain continuous-trace C*algebras and in particular certain subalgebras of interval algebras are inductive limits of special subalgebras.\\
\indent Let $A$ be a separable continuous-trace C*-algebra with spectrum $[0,1]$. In particular $A$ can be a full hereditary subalgebra of $M_n \otimes C[0,1]$. Clearly $H^{3}([0,1],\mathbb{Z})=H^{3}(\hat{A},\mathbb{Z})=0$. As a consequence of the Dixmier-Douady classification (cf. \cite{dix0}) we get that $A$ is stably isomorphic to $C[0,1]$. Moreover there exists $a:[0,1]\rightarrow \mathbb{K}$ such that $a\in A$ and $a(t)$ is a rank one projection for all $t$ in $[0,1]$. We notice that $A$ can be realized as a hereditary sub-C*-algebra of $\mathbb{K} \otimes C[0,1]$.\\
\indent We begin with two lemmas:
\begin{lemma} Let $p:[0,1] \rightarrow M_n$ be a lower semicontinuous projection-valued function that if composed with the rank function gives rise to the dimension function of a hereditary subalgebra $A$ of $M_n \otimes C[0,1]$. 
 Assume $p$ has constant rank equal $k$ on some open interval $(a,b)\subset [0,1]$ and $k \geq 1$. Then there exist a lower semicontinuous projection-valued function  $q:[0,1] \rightarrow M_n$ such that:\\
\indent (1)  $q(t)$ is a projection of rank $1$ for all $t \in (a,b)$ and zero otherwise;\\
\indent (2)  $p(t) \geq q(t)$ for all $t \in (a,b)$. 
\end{lemma}
\begin{proof}
  We know that $p(t)$ has constant rank equal $k$ on $(a,b)$. By Lemma 8.0.8 there exists a continuous family of unitaries $U:(a,b)\rightarrow M_n(\mathbb{C})$ such that:
 $$U(t)p(t)U(t)^{*}=P(t)=P,\;\;\mathrm{for}\;\;t\in(a,b).$$
\indent Because $A$ is a continuous trace C*-algebra with spectrum $[0,1]$,
 there exists $a:[0,1]\rightarrow M_n$ such that $a\in A$ and 
$a(t)$ is a rank one projection for all $t$ in $[0,1]$. Therefore for any $t \in (a,b)$ we have rank$(a(t))=1 \leq$ rank$(P(t))$. Thus there exists a family of unitaries $V:(a,b)\rightarrow M_n(\mathbb{C})$ such that $V(t)a(t)V(t)^{*} \leq P(t)=P$.\\
\indent Define $q(t)$ to be equal to $U(t)^{*}V(t)a(t)V(t)^{*}U(t)$ for $t \in (a,b)$ and equal to zero otherwise to get the desired conclusion of the lemma.
\end{proof}

\begin{lemma}
Let $p:[0,1] \rightarrow M_n$ be a lower semicontinuous projection-valued function that if composed with the rank function gives rise to the dimension function of a hereditary subalgebra $A$ of $M_n(\mathbb{C}) \otimes C[0,1]$. 
 Assume $p$ has maximal rank equal $k$ on some open interval $(a,b)\subset [0,1]$ and $k \geq 1$. Then there exist a lower semicontinuous projection-valued function \\ $q:[0,1] \rightarrow M_n(\mathbb{C})$ such that:\\
\indent (1)  $q(t)$ is a projection of rank $1$ for all $t \in (a,b)$ and zero otherwise;\\
\indent (2) $p(t)-q(t)$ is a lower semicontinuous projection-valued function.
\end{lemma}
\begin{proof} Because $p(t)$ is a lower semicontinuous projection-valued function and $a,b$ are points of discontinuity we have: $$\mathrm{rank}(p(a))<\mathrm{rank}(p(t))=k,\;\;t \in (a,b);$$
$$\mathrm{rank}(p(b))<\mathrm{rank}(p(t))=k,\;\;t \in (a,b).$$
\indent Next we choose two projections of rank $1$, $p_a$ and $p_b$, such that $p_a$ is orthogonal to $p(a)$ and $p_b$ is orthogonal to $p(b)$. Moreover $p_a$ and $p_b$ are homotopic equivalent. Therefore we can find a path of rank $1$ projections $a(t)_{t\in (a,b)}$ connecting $p_a$ and $p_b$. Now we proceed as in the previous lemma 8.1.1 and modify the path $a(t)$ to get $q(t)$ with properties that $q(t)\leq p(t)$ and $\mathrm{rank}(q(t))=1$ for all $t \in (a,b)$, zero otherwise. Furthermore $p(t)-q(t)$ is a lower semicontinuous projection-valued function as desired.  
\end{proof}

\begin{theorem} Let $p:[0,1]\rightarrow M_n(\mathbb{C})$ be a lower semicontinuous projection-valued function  which comes from the dimension function of a hereditary subalgebra $A$ of $C[0,1] \otimes M_n(\mathbb{C})$.
 Then $p$ is a sum of special lower semicontinuous projection-valued functions, where by special we mean that all projections in the range of the functions have rank at most $1$.
\end{theorem}
\begin{proof} Let $A_k:=\{t\in [0,1]|\;\mathrm{rank}(p(t))=k \}$ where $k$ is the maximum value of the rank for the projections in the range of $p$. Because $p$ is a lower semicontinuous function we get that $A_k$ is an open subset of $[0,1]$. Therefore $A_k$ is a countable disjoint union of $(I_i)_i$, where each $I_i$ is an open interval. Next we apply Lemma 8.1.1 to get a lower semicontinuous projection-valued function $p_{1}:[0,1] \rightarrow M_n(\mathbb{C})$ which is a rank $1$ projection on the union of $(I_i)_i$ and zero on the complement. In other words $p_1(t)$ is zero for $t$ not in $A_k$ but in $[0,1]$.

 By the lemma 8.1.1 $p(t)-p_1(t)$ is a lower semicontinuous projection-valued function. Moreover the maximum value of the rank of projections appearing in the range of $p(t)-p_1(t)$ is at most $k-1$.\\
\indent Next repeat the argument with $p(t)-p_1(t)$ instead of $p(t)$ to get another lower semicontinuous projection valued function $p_2(t)$ which has values projections of rank one precisely for the set $$A_{k-1}=\{t\in [0,1]:\mathrm{rank}(p(t)-p_1(t))=k-1\}.$$
\indent Because $p(t)$ has rank $k$ precisely on $A_{k}$ and $p_1$ is zero on the complement of $A_{k}$ we have $A_{k}\subset A_{k-1}$. Continuing in this way, after at most $k$ steps we obtain a decomposition for $p(t)$ as desired: $$p=\sum\limits_{i=1}^{k}p_i$$ 
 where each $p_i$ is zero except an open set $A_i$ where it takes values projections of rank $1$. Moreover for any $i\in \{1,\dots ,k-1 \}$ we have $A_{i+1} \subset A_{i}$.
\end{proof}

 We obtain a similar decomposition result for certain dimension functions:
\begin{prop} Any lower semicontinuous integer-valued function is a sum of lower semicontinuous functions whose values are either $0$ or $1$.\\
\indent This decomposition is unique if we ask for the partial sums to form an increasing sequence.
\end{prop} 
\begin{proof}We know that the dimension function (which assign to each point in the spectrum the dimension of that irreducible representation) is a lower semicontinuous positive, integer-valued function and bounded above by the highest dimensions of the irreducible representations of the given algebra. \\
\indent Notice that the maximum value of the dimension function is taken on an open subset of $[0,1]$. This is an important fact which is fundamental to obtain our decomposition. \\ 
\indent Consider the step function which takes the value 1 on the open set where the given lower semicontinuous map takes its maximum value. It is easy to see that this step function is lower semicontinuous and the difference between the dimension function and the just defined step function is lower semicontinuous. Moreover this new lower semicontinuous map has positive values and a smaller maximum value than the maximum of the previous function. Repeat this argument until we get a function with maximum $1$. We conclude that the original function is a sum of lower semicontinuous functions that are characteristic functions on some open subsets of $[0,1]$ (i.e. they take values either $0$ or $1$).\\
\indent Notice that the characteristic sets of the above sequence of step functions increases as we decrease the maximum value of the dimension function.\\ 
\end{proof}

\begin{remark} The dimension function of a full hereditary subalgebra takes values bigger than $1$. This in turn implies that the last function in the decomposition is $1$ everywhere on $[0,1]$.
\end{remark}
\indent A simple example of the decomposition given by the Proposition 8.1.4:

\begin{center}
\unitlength=1mm
\begin{picture}(110, 30)(-4, -4)
\put(-2, 7){\makebox(0, 0){1}}
\put(-2, 14){\makebox(0, 0){2}}
\put(17, -3){\makebox(0, 0){1}}
\put(25, 7){\makebox(0, 0){$=$}}
\put(-5, 0){\line(1,0){25}}
\put(0,7 ){\line(1,0){7}}
\put(0, -5){\line(0, 1){25}}
\put(7, 14){\line(1, 0){10}}
\put(17, -1){\line(0, 1){2}}
\put(38, 7){\makebox(0, 0){1}}
\put(57, -3){\makebox(0, 0){1}}
\put(66, 7){\makebox(0, 0){$+$}}
\put(35, 0){\line(1,0){25}}
\put(40, 7 ){\line(1,0){17}}
\put(40, -5){\line(0, 1){25}}
\put(7, 14){\line(1, 0){10}}
\put(57, -1){\line(0, 1){2}}

\put(72, 7){\makebox(0, 0){1}}
\put(91, -3){\makebox(0, 0){1}}
\put(70, 0){\line(1,0){25}}
\put(81,7 ){\line(1,0){10}}
\put(74, -5){\line(0, 1){25}}
\put(91, -1){\line(0, 1){2}}
\end{picture}
\end{center}
\vspace{4mm}

 The following lemma is used to prove the main result of this section:
\begin{lemma} Let $A$ be a C*-algebra and $a$ an element of $A$. Then the hereditary sub-C*-algebra generated by $a^{*}a$ is isomorphic to the hereditary sub-C*-algebra generated by $aa^{*}$.
\end{lemma}

\begin{proof} 
Denote with $H_1$ and $H_2$ the hereditary C*-subalgebra generated by $a^{*}a$ and respectively $aa^{*}$:
 $$H_1=\overline{a^{*}aAa^{*}a},\;\;H_2=\overline{aa^{*}Aaa^{*}}.$$

\indent Consider the polar decomposition $a=v(a^{*}a)^{1/2}=(aa^{*})^{1/2}v$ and observe that we also have:
  $$H_1=\overline{(a^{*}a)^{1/2}A(a^{*}a)^{1/2}},\;\;H_2=\overline{(aa^{*})^{1/2}A(aa^{*})^{1/2}}.$$

 Now define two maps: $$\phi :(aa^{*})^{1/2}A(aa^{*})^{1/2} \rightarrow (a^{*}a)^{1/2}A(a^{*}a)^{1/2}$$
            $$ \phi(b)=v^{*}bv,$$
            $$\psi :(a^{*}a)^{1/2}A(a^{*}a)^{1/2} \rightarrow  (aa^{*})^{1/2}A(aa^{*})^{1/2}$$
            $$ \psi(c)=vcv^{*}.$$
\indent We use $v^{*}(aa^{*})^{1/2}=(a^{*}a)^{1/2}v^{*}$, $(aa^{*})^{1/2}v=v(a^{*}a)^{1/2}$ to see that $\phi$ and $\psi$ are maps between the specified sets: $$\phi((aa^{*})^{1/2}m(aa^{*})^{1/2})=v^{*}(aa^{*})^{1/2}m(aa^{*})^{1/2}v=$$
$$=(a^{*}a)^{1/2}v^{*}mv(a^{*}a)^{1/2}\in (a^{*}a)^{1/2}A(a^{*}a)^{1/2},\;\;m\in A.$$
  And similarly for $\psi$.\\
\indent To show that $\phi$ and $\psi$ are *-homomorphisms and their composition is identity we use that $v^{*}v$ and $vv^{*}$ are projections on the ranges of $(a^{*}a)^{1/2}$ and respectively $(aa^{*})^{1/2}$:
 $$\phi(b_1b_2)=\phi((aa^{*})^{1/2}m_1(aa^{*})^{1/2}(aa^{*})^{1/2} m_2(aa^{*})^{1/2})=$$ 
$$=v^{*}(aa^{*})^{1/2}m_1(aa^{*})^{1/2}(aa^{*})^{1/2}m_2(aa^{*})^{1/2}v=$$
$$=v^{*}(aa^{*})^{1/2}m_1(aa^{*})^{1/2}vv^{*}(aa^{*})^{1/2}m_2(aa^{*})^{1/2}v=$$
         $$=\phi(b_1)\phi(b_2),$$
         $$\psi(\phi(b))=\psi(v^{*}bv)=vv^{*}bvv^{*}=b,$$
  $\mathrm{where}\;b_{i}=(aa^{*})^{1/2}m_{i}(aa^{*})^{1/2},\;i\in\{1,2\},\;m_i \in A.$  
 $$\phi(\psi(c))=\phi(vcv^{*})=v^{*}vcv^{*}v=c,$$
 with $c_{j}=a^{*}am_{j}a^{*}a,\;j \in \{1,2\},\;m_j\in B.$\\
\indent Finally we extend the isomorphism to the closure of the subalgebras to obtain the desired isomorphism.

\end{proof}

\begin{theorem} Let $A$ and $B$ be two full hereditary sub-C*-algebras of $C[0,1]\otimes \mathbb{K}$, with the same dimension function $f$, where we also keep track of the identification of the spectrum of the algebras $A$ and $B$. Then there is an isomorphism between $A$ and $B$ which preserves the identification of the spectrum.    
\end{theorem}

\begin{proof}
Let $P^{A}$ and $P^{B}$ be the two lower semicontinuous projection-valued functions that gives rise to the dimension function $f$ (see Remark 7.0.21). Because the dimension function is the same we have $$\mathrm{rank}(P^{A}(t))=\mathrm{rank}(P^{B}(t))\;\;\mathrm{for}\;\mathrm{all}\;t\in [0,1].$$
\indent The idea of the the proof is to construct a partial isometry $V$ such that $V^{*}V=P^{A}$ and $VV^{*}=P^{B}$ and to apply the Lemma 8.1.6.\\  
\indent We begin by applying Theorem 8.1.3 to get a decomposition of the lower semicontinuous functions  $$P^{A}=\sum\limits_{i=1}^{N}P_{i}^{A},$$
           $$P^{B}=\sum\limits_{i=1}^{N}P_{i}^{B}.$$
 Note that $N$ is the maximum value of the dimension function. Moreover for each $i\in \{1, \dots, N \}$, $P^{A}_i$ and $P^{B}_i$ are zero except the same subset $S_i$ of $[0,1]$. For $t\in S_i$ we have $$\mathrm{rank}(P^{A}_{i}(t))=\mathrm{rank}(P^{B}_{i}(t))=1.$$ 
 As shown by Glimm in \cite{gli} for each $t_0 \in [0,1]$ we can find $W_{i}(t)\in C[0,1]\otimes \mathbb{K}$ such that $W_{i}(t_{0})W^{*}(t_{0})_i=P^{A}_i(t_{0})$, $W_i^{*}(t_{0})W(t_{0})_i=P^{B}_i(t_{0})$ and $||W_{i}(t)P^{A}(t)W_{i}(t)^{*}-P^{B}(t)||< \frac{1}{2}$ for all $t$ in a neighborhood of $t_0$.\\
\indent  Making use of a continuous function $g:[0,1]\rightarrow [0,2]$ defined to be zero in a neighborhood of zero and $\frac{1}{t}$ for $t\geq \frac{1}{2}$ we can define a continuous family of partial isometries $$\widehat{W_{i}}(t)=P^{A}(t)[g(P^{A}(t)W_{i}(t)P^{B}(t)W_{i}(t)^{*}P^{A}(t)]^{1/2}W_{i}(t)P^{B}(t)W_{i}(t)^{*}$$ such that for all $t$ in a neighborhood of $t_0$ we have: $$\widehat{W_{i}}(t) \widehat{W_{i}}^{*}(t)=P^{A}_{i}(t)$$ $$\widehat{W_{i}}^{*}(t) \widehat{W_{i}}(t)=P^{B}_{i}(t).$$ 
\indent As in the Lemma 8.0.4 we extend the partial isometry on each subinterval of the set $S_i$. Therefore for each open interval $I \subseteq S_i$ we find a family of partial isometries $(V_i(t))_{t\in I}$ such that $V_i(t)V^{*}_{i}(t)=P^{A}_i(t)$, $V_i^{*}(t)V_{i}(t)=P^{B}_i(t)$. For $t \in [0,1]$ but not in $S_i$ we let $V_i(t)=0$. Notice that for any $i \in \{1, \dots ,N\}$ we have: $$V_i(t)V^{*}_{i}(t)=P^{A}_i(t),\;\forall t \in [0,1],$$ $$V_i^{*}(t)V_{i}(t)=P^{B}_i(t),\;\forall t \in [0,1].$$

 Let $t\in [0,1]$. Define $$V(t)=\sum\limits_{i=1}^{N}V_i(t).$$
\indent In this way we get a family of partial isometries $V(t)_{t\in[0,1]}$ such that $V(t)V(t)^{*}=P^{B}(t)$ and $V(t)^{*}V(t)=P^{A}(t)$ for all $t \in [0,1]$. We notice that the continuity of a function $a(t)$, such that $P^{A}(t)a(t)=a(t)$ implies the continuity of $V(t)a(t)V(t)^{*}$. It is a continuous function because $V_ia(t)$ is continuous for any $i \in \{1, \dots, N\}$. Moreover we have the following polar decomposition: $V(t)= V(t)(V(t)^{*}V(t))^{1/2}$.\\
\indent Now we apply Lemma 8.1.6 with $V(t)$ instead of $a$ to obtain the desired isomorphism: $$A \cong B.$$
\end{proof}

\begin{corollary} Let $A$ and $B$ be two continuous-trace sub-C*-algebras of $C[0,1]\otimes M_n$, with the same spectrum the set $[0,1]$ and the same finite step function dimension function, where we also keep track of the identification of the spectrum of the interval algebra. Then there is an isomorphism between $A$ and $B$.
\end{corollary}
\begin{proof} The idea of the proof is to reduce the problem to the case of full hereditary subalgebras.\\
\indent Since $H^{3}([0,1],Z)=0$ and $[0,1]$ is the spectrum of continuous-trace sub-C*-algebras $A$ and $B$, by the Dixmier-Douady classification (\cite {dix0}), we get that both $A$ and $B$ are stably isomorphic to $C[0,1]$. Therefore $A$ and $B$ can be realized as full hereditary subalgebras of $C[0,1]\otimes \mathbb{K}$.\\
\indent To get the isomorphism we apply Theorem 8.1.7.
\end{proof}

\indent  As a consequence of the Theorem 8.1.7 we prove the following approximation result:

\begin{theorem} Any hereditary sub-C*-algebra of an interval algebra can be realized as the closure of an increasing union of step hereditary sub-C*-algebras. 
\end{theorem}
\begin{proof}
 By Theorem 8.1.7 there is a one to one correspondence between hereditary subalgebras (considered up to isomorphism) and lower semicontinuous projection-valued functions. Therefore to prove the result it is enough to approximate a general lower semicontinuous projection-valued function with an increasing sequence of finite step function lower semicontinuous projection-valued functions.\\
\indent Let $f:[0,1]\rightarrow M_n$ be a lower semicontinuous projection-valued function. By the Theorem 8.1.3 we get that $$f= \sum\limits_{i=1}^N f_i$$
 where each $f_i$ is a lower semicontinuous projection-valued function, with all projections in the range of rank at most $1$.\\
\indent Recall that any open subset of $[0,1]$ can be approximated with the increasing union of finitely many open intervals. Therefore each $f_i$ can be approximated with finite step function projection-valued function: $$f_i= \lim_{n} f_{i}^{n}$$
 Then $$f = \lim_{n} \sum\limits_{i=1}^{N} f_{i}^{n},$$
where each $$\sum\limits_{i=1}^{N} f_{i}^{n}$$ is a finite step function projection-valued function.\\
\indent The hereditary sub-C*-algebras generated by $$\sum\limits_{i=1}^{N} f_{i}^{n}$$ define the desired approximation.
\end{proof}
\begin{corollary} Any continuous-trace C*-algebras whose spectrum is $[0,1]$ can be realized as the closure of an increasing union of special continuous-trace C*-algebras, special in the sense that their dimension function is a finite step function.
\end{corollary}

\section{ C*-algebras with a finite presentation and stable relations}

 It is important to know that certain algebras are finitely presented and have stable relations because we can modify maps which are close to be *-homomorphisms into exact *-homomorphisms. \\
 \indent We are particularly interested in so called weakly stable relations. Roughly speaking, weak stability means that close to an approximate representation there is an exact representation, in any C*-algebra. Formally, we have the following definition given by Loring in \cite{lor0}: 
\begin{definition} Suppose $A$ is a finitely presented C*-algebra: $A=C^{*}<G,R>$. If, for every $\delta > 0$, there exists $\epsilon >0$ and 
 $$\sigma_{\epsilon}: C^{*}<G,R>\rightarrow C^{*}_{\epsilon}<G,R>$$ such that $$||\sigma_{\epsilon}\circ P_{\epsilon}(g^{\epsilon}_{j})-g^{\epsilon}_{j}||\leq \delta $$ then the set of relations $R$ is called weakly stable.
\end{definition}

\begin{remark} If, in addition in the above definition, $P_{\epsilon}\circ \sigma_{\epsilon} = \mathrm{id}$, via a homotopy $\psi$, such that $$||\psi_{t}(g_j)-g_j|| \leq \delta $$ then $R$ is called strongly stable. If in fact $P_{\epsilon}\circ \sigma_{\epsilon} = \mathrm{id}$, then $R$ is exactly stable.
\end{remark}

\section{Special continuous-trace C*-algebras are finitely presented and have stable relations}
 
\indent Our goal in this section is to present a procedure to obtain a complete finite set of generators and a complete finite set of relations for special sub-C*-algebras of interval algebras.\\
\indent Let $A$ be a full step-hereditary subalgebra of $M_n \otimes C[0,1]$ with $f$ its finite step dimension function.\\
\indent The idea is to construct another concrete hereditary subalgebra of $M_n \otimes C[0,1]$, isomorphic to $A$, the advantage being that for this newly constructed algebra it is easy to write a finite set of generators.\\
\indent  By Proposition 8.1.1 the dimension function has the decomposition: $$f = \sum\limits_{i=1}^{N} f_i$$ where $f_i$ is the characteristic function of some open set $A_i$ of $[0,1]$ and each $A_i$ is a finite disjoint union of open intervals. We also have $A_i \subseteq A_{i+1}$.\\ 
\indent Next we construct a concrete hereditary sub-C*-algebra of $M_n \otimes C[0,1]$ with the dimension function equal to $f$. \\
\indent The hereditary sub-C*-algebra $H$ generated by the block diagonal form:

 $$ \left( \begin {array} {cccc} C_{0}(A_1) &     & 0 \\
                                     & C_{0}(A_2) &   \\
                                     &     & \ddots  \\
                                  0  &     &  C[0,1]  \\
 \end{array} \\ \right) \subseteq M_n \otimes C[0,1]  $$
has the dimension function $f = \sum\limits_{i=1}^{N}f_i$.\\
\indent We notice that all continuous-trace C*-algebras with spectra the closed interval $[0,1]$ and all irreducible representations finite dimensional are isomorphic to one of the very special form namely diagonal generated form and that for any dimension function we can construct such a model.\\
\indent Recall that a hereditary sub-C*-algebra $H$ of some C*-algebra $M$ has the property:  whenever $x \in M,\;\;xx^* \in H\;\;\mathrm{and}\;\; x^{*}x \in H \;\;\;\mathrm{then}\;\;\;x \in H.$\\
\indent Therefore we get that the concrete hereditary sub-C*-algebra $H$ contains the following sub-C*-algebra $B$:
 $$ \left( \begin {array} {cccccc} C_{0}(A_1) & C_{0}(A_1) & C_{0}(A_1) & \dots &C_{0}(A_1)\\
                                  C_{0}(A_1)& C_{0}(A_2) & C_{0}(A_2)& \dots & C_{0}(A_2)\\
                                  C_{0}(A_1)& C_{0}(A_2)    &  C_{0}(A_3) & \dots & C_{0}(A_3) \\
                                   \vdots   &  \vdots & \vdots  & \ddots     \\
                                 C_{0}(A_1)& C_{0}(A_2)& C_{0}(A_3)& \dots   & C[0,1]  \\
 \end{array} \\ \right) \subseteq M_n \otimes C[0,1].  $$
\vspace{4mm}
\indent We notice that subalgebra $B$ is a hereditary subalgebra of $M_n \otimes C[0,1]$: $$b_1ab_2 \in B\;\;\mathrm{whenever}\;b_1,b_2 \in B,\;a\in M_n\otimes C[0,1].$$
\indent Therefore $B=H$.\\
\indent  Moreover since the dimension function of $H$ is $f$, by Theorem 8.1.7 we get: $$A \cong H .$$
\subsection { Finding a finite set of generators for the concrete algebra}
\indent We consider the subset of elements with nonzero entries on the last row:
 $$ \left( \begin {array} {cccccc} 0 & 0 & 0 &\dots & 0  \\
                                  0 & 0 & 0 & \dots & 0\\
                                  0& 0    &  0 & \dots & 0 \\
                                   \vdots   &  \vdots & \vdots  & \ddots     \\
                                 C_{0}(A_1)& C_{0}(A_2)& C_{0}(A_3)& \dots   & C[0,1]  \\
 \end{array} \\ \right) \subseteq M_n \otimes C[0,1].  $$ 
\vspace{2mm}
\indent For the first block entry $C_{0}(A_1)$ we consider the following matrix valued functions:
 $$f_{n1}^{1}\otimes e_{n1},f_{n1}^{2}\otimes e_{n1},\dots, f_{n1}^{c_{1}}\otimes e_{n1}$$
  where by $c_{1}$ we denote the number of connected components of the open set $A_1$, $f_{n1}^{i}$ are non-zero functions on the corresponding open intervals of $A_1$ and $e_{n1}$ is a matrix unit.\\
\indent For instance if an open interval is $(a,b)$ we consider the function $t \rightarrow |e^{2\pi i \frac{t-a}{b-a}}-1|$ defined on $(a,b)$. Notice that this function is zero at the end-points of the support set.\\
\indent We may have half-open intervals at the end points of the intervals: $[0,a)$ or $(b,1]$ in which case we consider the function $f(t)=a-t$ respectively $f(t)=t-b$.\\
 \indent For the second block $C_{0}(A_2)$ on the last row we get:
 $$f_{n2}^{1}\otimes e_{n2},f_{n2}^{2}\otimes e_{n2},\dots, f_{n2}^{c_{2}}\otimes e_{n2}$$
 where $f_{n2}^{i}$ are non-zero functions on the corresponding open intervals of $A_2$ and $e_{n2}$ is a matrix unit. As before $c_{2}$ is the number of connected components of $A_2$. \\
\indent Continuing this procedure we obtain $$f_{ni}^{1}\otimes e_{ni},f_{ni}^{2}\otimes e_{ni},\dots, f_{ni}^{c_{i}}\otimes e_{ni},\;\mathrm{for}\;i\in \{1,\dots ,n-1\},$$
 where by $c_i$ we denote the number of connected components of the open sets $A_i$ appearing in the decomposition.

\indent For the last block entry $C[0,1]$ we consider: $$f_{nn}\otimes e_{nn}, {\bf 1}\otimes e_{nn}.$$
 \indent Note that $f_{nn}^{c_n}=t$, $c_{n}=1$ and ${\bf 1}(t)=1$ for all $t\in [0,1]$.\\
\indent We observe that $({\bf 1}) \otimes e_{nn}$ is a full abelian projection.\\
\indent This proposed set of generators is finite because the sets $A_i$ are finite union of open intervals.\\
\indent In order to see that the proposed set of elements generates the concrete algebra we use the Kaplansky's Stone-Weierstrass Theorem (see \cite {dix0}). In other words we prove that the C*-algebra generated by the proposed set of elements separates points of the set of all pure states of $A$ together with $0$. \\
\indent For instance let $\phi$ and $\psi$ be two pure states of the C*-algebra $A$. Now irreducible representations of $A$ are evaluations at points in the interval $[0,1]$. Assume that the pure state $\phi$ corresponds to the irreducible representation $(\pi _{\phi},H_{\phi})$ which is an evaluation at the point $t_{\phi} \in [0,1]$ and the pure state $\psi$ corresponds to the irreducible representation $(\pi _{\psi},H_{\psi})$ which is an evaluation at the point $t_{\psi} \in [0,1]$. If $dimH_{\phi}>dimH_{\psi}$ then the two states are separated by an element of the proposed set of generators which at $t_{\phi}$ is not zero. For the case when $dimH_{\phi}=dimH_{\psi}$ we observe that $t(f^{i}_{nj}\otimes e_{nj})$ belongs to the C*-algebra generated by the proposed set of generators. This observation can be verified either by using the Dauns-Hofmann theorem or by doing some algebra computations. As a consequence of the above observation we have that if $\phi$ and $\psi$ are not separated by the proposed set of generators then they are separated by some $t(f^{i}_{nj}\otimes e_{nj})$. Moreover this also shows the separation between $0$ and any other pure state of $A$.

\subsection { Finding a finite set of abstract generators and relations}
 \indent We show that a step-hereditary sub-C*-algebra $A$ can be described as a universal C*-algebra generated by a finite set of generators that satisfy a finite set of relations. We notice that there are relations which are neither norm conditions nor *-polynomials relations.\\   
\indent In the previous paragraph we found a concrete description of the algebra $A$ in terms of the boundary conditions, the advantage being that we recognize an explicit finite set of concrete generators. Thus the cardinal of the set of ``abstract generators'' is known.\\
\indent By what we said above we get the following set of abstract generators:
 $$G=\{ P,\;X_{i}^{k},\;\;i \in \{1,\dots n \},n\geq 2\;\;k \in \{ 1 \dots c_i \}\}.$$ 
\indent Assume that $A^{k}_{i}$ are open subsets of $[0,1]$, $i \in \{1,\dots n\},\;k \in \{1,\dots c_{i}\}$ such that $c_n=1$, $A_n=[0,1]$ and $\bigcup\limits_1^{c_i}A^{k}_{i} \subseteq \bigcup\limits_{1}^{c_{i+1}}A^{k}_{i+1}$ for all $i \in \{1,\dots n\}$.\\
\indent By working with the concrete generators we find the following groups of relations.\\
\indent  A first group of relations is:  $$(R_1):\;\; P=P^*=P^2,\;\;||X_{i}^{k}||\leq 2$$
 $$||X_{n}^{c_n}|| \leq 1.$$
\indent A second group of relations is:  $$(R_2):\;\; X_{i}^{k}X_{j}^{l}=0,\;,\; i,j \in {\{1,\dots n-1 \}},  k\in {\{1,\dots c_i \}, l \in \{1,\dots c_j \} }$$ 
 $$ X_{i}^{k}(X_{j}^{l})^{*}=0,\;i \neq j, \; i,j \in {\{1,\dots n-1 \}},  k\in {\{1,\dots c_i \}, l \in \{1,\dots c_j \} }$$ 
 $$X_{i}^{k}(X_{i}^{l})^{*}=0,\; i\in {\{1,\dots n-1 \}},\;k\neq l,\;\;  k,l\in {\{1,\dots c_i \} }$$ 
 $$PX_{i}^{k}=X_{i}^{k}\;\;\mathrm{and}\;X_{i}^{k}P=0,\;i \in {\{1,\dots n-1\}},\;k\in {\{1,\dots c_i \}}$$ 
 $$PX_{n}=X_{n}P=X_{n}$$
$$0 \leq X_n \leq 1.$$

\begin{prop} There is an isomorphism
 $$C^{*}<G|R_1 \cup R_2 \cup R_3> \cong A$$
 sending 
 $$X_{i}^{k} \mapsto f_{i}^{k}(t) \otimes e_{ni},\;\;i \in {\{1,\dots n\}},k\in {\{1,\dots c_i \}}$$ 
 $$P \mapsto {\bf 1} \otimes e_{nn}.$$
\end{prop}

\begin{proof} We shall construct the universal C*-algebra in the following way (cf. \cite{lor}). Denote with $R_{||*||}$ the norm relations $\{||X_1|| \leq 1, ||X^{i}_{k} \leq 2\}$; with $R_{p}$ the set of *-polynomials relations; and with $R_{f}$ the set of relations of the form $f(x)=y$. Then $$C^{*}<G|R_{||*||}>=F_{n}$$ where $n$ is the number of generators and $$C^{*}<G|R_{||*||}\cup R_{p}>=F_{n}/I$$ where I is the ideal generated by the *-polynomials in $R_{p}$. Now in $F_{n}/I$ the elements of $R_{f}$ all make sense, so we let $J$ equal the ideal generated by $R_{f}$ and set $$C^{*}<G|R_1 \cup R_2 \cup R_3>=C^{*}<G|R_{||*||} \cup R_{p} \cup R_{f}>=(F_N/I)/J.$$
\indent The universality condition implies that there exists a surjective map 
$$C^{*}<G|R_1 \cup R_2 \cup R_3> \rightarrow A$$
 that we show is an isomorphism.\\
\indent Since $PX_{i}^{k}=X_{i}^{k},\mathrm{for}\;\mathrm{all}\; i \in {\{1,\dots n\}, k \in {\{1,\dots c_i \}}} $ it follows that $P$ is a full projection.\\
\indent Making use of relations $(R_2)$ we obtain that $P$ is an abelian projection:
 $$PX_{i_{1}}^{k_{1}}X_{i_{2}}^{k_{2}} \dots X_{i_{m}}^{k_{m}}P=0$$ if there are no adjoints of the generators appearing in the product. Moreover to get nonzero such products of generators we have to have:
  $$PX_{i_{1}}^{k_{1}}(X_{i_{1}}^{k_{1}})^{*} X_{i_{2}}^{k_{2}}(X_{i_{2}}^{k_{2}})^{*}\dots X_{i_{m}}^{k_{m}}(X_{i_{m}}^{k_{m}})^{*}P=$$
  $$PX_{i_{1}}^{k_{1}}(X_{i_{1}}^{k_{1}})^{*}PP X_{i_{2}}^{k_{2}}(X_{i_{2}}^{k_{2}})^{*}P\dots PX_{i_{m}}^{k_{m}}(X_{i_{m}}^{k_{m}})^{*}P=$$
  $$f_{i_{1}}^{k_{1}}(X_n)f_{i_{2}}^{k_{2}}(X_n)\dots f_{i_{m}}^{k_{m}}(X_n) \in C^{*}(X_n,P) \cong C[0,1]$$
\indent By the above relations we have that the spectrum of the universal C*-algebra is a subset of the set $[0,1]$.\\
\indent In order to have the desired isomorphism we still have to prove that for any $\alpha \in [0,1]$ there is a corresponding irreducible representation $\pi _{\alpha}$ of the set $G$. We construct $\pi _{\alpha}$ as follows
 $$P \mapsto {\bf 1} \otimes e_{nn},$$ 
$$X_{i}^{k} \mapsto \sqrt{f_{i}^{k}(\alpha)} \otimes e_{ni},\;\;i \in {\{1,\dots n\}},k\in {\{1,\dots c_i \}}$$ 
and the Hilbert space $H_{\alpha}$ is $\mathbb{C}^{n}$. In checking that $\pi_{\alpha}$ is irreducible we distinguish the following cases: if $\alpha \in A^{k}_1$ then $dimH_{\alpha}=n$, where by $A^{k}_1$ is the k-th connected component of $A_1$. Notice that the closure of $A^{k}_1$ is the support of $f^{k}_{n1}$. In the case that $\alpha \in A^{k}_{i_{0}}$ but $\alpha \notin A^{k}_{i_{0}-1}$ we obtain an irreducible representation of dimension $n-i_{0}+1$ because $\sqrt{f^{k}_{ni}(\alpha)}=0$ which implies $X_{i}^{k}$ are sent to $0$ for $i<i_0$.
\end{proof}
\begin{remark} The above proof shows that the norm relations $||X^{k}_{i}||\leq2$ are redundant except the norm relation $||X_n||\leq 1$. 
\end{remark}

 Now that we have identified the universal C*-algebra for the set of generators $G$ and relations $R_1 \cup R_2 \cup R_3$, the next step is to prove that these relations are stable relations. We start with two lemmas:

\begin{lemma} For every $\epsilon$, there is a $\delta$ such that, given operators $x_1,x_2,\dots,x_n,p$ in a C*-algebra $\textit A$ with relations $R_1$ and $R_2$ satisfied within $\delta$, there exist operators $\widehat{x_1},\widehat{x_2},\dots,\widehat{x_n},\hat{p}$ in $\textit A$, which exactly satisfy the relations $R_1$, $R_2$, and $||x_i-\widehat{x_i}||\leq \epsilon,\;i\in {\{1, \dots ,n\}},\;||p-\hat{p}||\leq \epsilon $.
\end{lemma}

\begin{proof} First we perturb $p$ like in \cite{lor} (Lemma 4.2.2) to obtain a projection $\hat{p}$ close to $p$.\\
\indent Next define $\widehat{x_i}=\hat{p}x_i(1-\hat{p}),\;i \in {\{1, \dots ,n-1\}}$ and $\widehat{x_n}=\hat{p}x_n\hat{p}$. It is easy to see that all relations $R_1$ and $R_2$ are satisfied and:
 $$||x_i-\widehat{x_i}||=||x_{i}-\hat{p}x_{i}(1-\hat{p})||= $$
 $$=||x_i -px_{i}(1-p)+px_{i}(1-p)- \hat{p}x_{i}(1-p)+ \hat{p}x_{i}(1-p)-\hat{p}x_{i}(1-\hat{p})|| \leq $$
$$\leq ||x_i -px_{i}(1-p)||+||\hat{p}-p||\;||x_i||\;||1-p||+$$
$$+ ||p-\hat{p}||\;||x_i||\;||(1-p)-(1-\hat{p})||\leq$$
$$\leq ||x_i-px_i||+||px_ip||+ 2\delta  +2 (\delta)^{2} \leq $$
$$\leq \delta + \delta + 2 \delta +2 (\delta)^{2} \leq \epsilon ,$$
 for suitable $\delta$. Also 
 $$||x_n-\widehat{x_n}||= ||x_{n} - pX_{n}p + px_{n}p -\hat{p}x_{n}p +\hat{p}x_{n}p -\hat{p}x\hat{p}||\leq$$
 $$\leq ||x_{n}-px_n||+||p||\;||x_n-x_{n}p||+||p-\hat{p}||\;||x_n||\;||p||+||\hat{p}||\;||x_{n}||\;||p-\hat{p}|| \leq $$
 $$ \leq \delta + \delta + 2 \delta + 2\delta \leq \epsilon.$$
 By choosing $\delta$ sufficiently small we get that all the perturbations $\widehat{x_i}$ are as small as needed.
\end{proof}

\begin{lemma} The non-polynomial relation $xx^*=f(y)$, where f is a continuous positive function on $[0,1]$, together with $||x||\leq 2$ and $0\leq y \leq 1$ is a stable relation:\\
\indent For all $ \epsilon \geq 0$ there exists $\delta \geq 0$ such that given two elements $x$ and $y$ in a C*-algebra $A$ with  $||xx^*-f(y)||\leq \delta $, $||x||\leq 2$ and $0\leq y \leq 1$,
there exists $\hat{x}$ and $\hat{y}$ in the C*-algebra $A$ such that:
 $$||x-\hat{x}||\leq \epsilon,\;||y-\hat{y}||\leq \epsilon \;\mathrm{and}$$
 $$\hat{x}(\hat{x})^{*}=f(\hat{y})$$
\end{lemma}

\begin{proof}
 
\indent Let $(a_n)$ be an approximate unit of the hereditary C*-algebra generated by $f(y)$. Then the following equality holds: $\lim\limits_{\rightarrow}a_{n}f(y)a_{n}=f(y)$. Thus there is some $N_{\epsilon}$ such that for $n \geq N_{\epsilon}$ we have $$|a_nf(y)a_n-f(y)| < \delta.$$
\indent Now we can change slightly $x$ so that $xx^{*}$ belongs to the hereditary C*-algebra generated by $f(y)$: $$\hat{x}=a_nx,\;n>N_{\epsilon}.$$
\indent Clearly $(\hat{x}\hat{x})^{*}$ belongs to the hereditary algebra generated by $f(y)$ and $$||(\hat{x}\hat{x})^{*}-f(y)||=||a_nxx^*a_n-f(y)||\leq $$
$$\leq ||a_nxx^*a_n-a_nf(y)a_n||+||a_nf(y)a_n-f(y)||\leq \delta +\delta $$
\indent Choose a spectral projection $q$ of $f(y)$ such that $a_n q=q$ for sufficiently large $n$. \\
\indent Since $xx^{*}$ is close to $f(y)$ we have: $$xx^{*}\geq f(y)-\delta .$$
\indent It follows that $$a_nxx^{*}a_n \geq a_nf(y)a_n-\delta (a_{n})^{2}.$$
\indent On the other hand there exists $M>0$ such that: $$a_nf(y)a_n \geq M (a_{n})^{2}.$$ 
\indent We conclude that $$(\hat{x}\hat{x})^{*}=a_nxx^{*}a_n \geq (M-\delta)(a_{n})^{2}\geq (M-\delta)(a_{n})^{2}q=(M-\delta)q.$$
\indent Hence we can find a continuous function $g$ which is not zero on the support of the spectral projection $q$, $g(y)$ is close to $y$ and moreover $f(g(y))$ is an element of the hereditary C*-algebra generated by $(\hat{x}\hat{x})^{*}$.\\
\indent  Let $\hat{y}=g(y)$ and $\widehat{\hat{x}}=vf(\hat{y})^{1/2}$ where $v$ is a partial isometry that appears in the polar decomposition of $\hat{x}$, $\hat{x}=v((\hat{x}\hat{x})^{*})^{1/2}$.\\
\indent Then $\widehat{\hat{x}}$ and $\hat{y}$ satisfy all the desired conditions.
\end{proof}
\indent We now prove the main stability result of this section:
\begin{theorem} Special continuous-trace C*-algebras have stable relations.
\end{theorem}
\begin{proof} The proof is an application of the previous lemmas: first apply Lemma 9.2.3 to fix the relations $R_1\cup R_2$. Next apply Lemma 9.2.4 for each relation from $R_3$. Note that $x_n$ in Lemma 9.2.3 is thought as being $y$ in Lemma 9.2.4. After we apply Lemma 9.2.4 we still have the relations from Lemma 9.2.3 exactly satisfied because $\hat{p}$ commutes with $y$.
\end{proof}
\begin{corollary} Full step hereditary sub-C*-algebras of interval algebras have stable relations.
\end{corollary}
 
\section{Inductive limits of special algebras}
\indent Now we are ready to show that the inductive limit of continuous-trace C*-algebras with spectra $[0,1]$ is isomorphic to the inductive limit of special continuous-trace C*-algebras. \\
\indent We start by proving a more abstract result. Denote with $C$ the class of algebras which have finite presentation and stable relations. Next consider two classes: the first one, $C'$, is the class of inductive limits of sequences of algebras from $C$ and the second one, $C''$ is the smallest class of algebras containing $C$ and is closed under taking inductive limits. Then we have the following theorem:
\begin{theorem} The class of algebras $C'$ and the class of algebras $C''$ are equal:
                       $$C'=C''$$ 
\end{theorem}

\begin{proof}  First we notice that one inclusion is obvious, namely $C' \subset C''$. The other inclusion is an application of Loring's theorem. We apply the theorem several times, depending on how many inductive limits operations we need to consider in order to obtain a given algebra in $C''$. For clarity we reproduce the Loring's theorem below:\\
\indent({\bf \cite{lor0}, Theorem 3.8}): Suppose $A$ is a C*-algebra containing a (not necessarily nested) sequence of sub-C*-algebras $A_n$ with the property: for all $\epsilon > 0$ and for any finite number of elements $x_1, \dots,x_k$ of $A$, there exist an integer $n$ such that:$$\{ x_1, \dots, x_k \} \subset_{\epsilon}A_n.$$ 
\indent If each $A_n= C^{*}<G_n,R_n>$, for all $n$, and $R_n$ is weakly stable, then $$A \cong \lim\limits_{\rightarrow}(A_{n_k}, \gamma_{k})$$
for some subsequence, and some connecting maps $\gamma_k:A_{n_k}\rightarrow A_{n_{k+1}}$.
 To see how the maps $\gamma_k:A_{n_k}\rightarrow A_{n_{k+1}}$ are defined we notice, as it was shown in \cite{lor}, that whenever we have a *-homomorphism $\gamma:A \rightarrow B$ and $C$ a sub-C*-algebra of $B$ such that $C$ almost contains the image under $\gamma$ of a finite subset of $A$, then there exists a *-homomorphism $\widehat{\gamma}:A \rightarrow C$ which is close to the original *-homomorphism $\gamma:A \rightarrow B$. Here ``almost'' and ``close'' should be interpreted in a suitable sense which is made precise in \cite{lor}.    
\end{proof}

\indent As a particular case of the above theorem we obtain:
\begin{theorem} Let $A$ be a simple inductive limit of continuous-trace C*-algebras whose building blocks have their spectrum homeomorphic to $[0,1]$. Then $A$ is an inductive limit of special continuous-trace C*-algebras whose building blocks have their spectrum homeomorphic to $[0,1]$ (in particular inductive limit of step hereditary sub-C*-algebras of interval algebras).
\end{theorem}
\begin{proof}
\indent Let $A=\lim\limits_{\rightarrow}A^{k}$, where $A^{k}$ are continuous-trace C*-algebras whose spectrum is an interval or a finite union of intervals (or hereditary sub-C*-algebras of interval algebras).\\
\indent We approximate each $A^{k}$ with special continuous-trace C*-algebras, namely continuous-trace C*-algebras whose dimension function is a finite step function. This can be done using Theorem 8.1.9.\\
\indent Hence we have

 $$\begin{array}{cccccccccc}

  \dots &\rightarrow & A^{k} &\rightarrow &A^{k+1} &\rightarrow&\dots &\lim\limits_{\rightarrow}A^{*}& \cong & A \\
          &  & \wr & & \wr & & & & \\
            & & \cup A^{k}_{i}& & \cup A^{k+1}_{i}& & \dots &  \\
          & & \uparrow & & \uparrow & & \\
            &  & \vdots & &\vdots &   &  \\
           & & \uparrow & & \uparrow & & \\ 
            & &  A^{k}_{i}& \rightarrow & A^{k+1}_{i} & & &    \\
             & & \uparrow & & \uparrow & & \\
 
  \end{array}$$
\indent In section 10 (Theorem 10.2.4) we show that special continuous-trace C*-algebras have finite presentation and stable relations. By Loring's lemma we conclude that:
  $$A \cong \lim\limits_{\rightarrow} A^{k}_{i} $$
\indent This last result can be viewed as a diagonal procedure to realize the algebra $A$ as an inductive limit of special algebras.\\
\indent  Therefore $A$ is an inductive limit of special continuous-trace C*-algebras (or inductive limit of step hereditary sub-C*-algebras of interval algebras).
\end{proof}
 
\section{The isomorphism theorem}
\indent In this section we prove the main classification result. In order to do this we use results from the previous section and results appearing in \cite{ste0}.\\

\indent Let $A$ a simple inductive limit of continuous-trace C*-algebras whose spectrum is $[0,1]$. Denote by $\Gamma(A)$ the invariant associated to $A$:
$$\Gamma(A)=[D(A), AffT^{+}A,Aff'A]$$

\indent Then Theorem 2.0.1 can be restated as follows:\\

\textit{ If $A$ and $B$ are simple inductive limits of continuous-trace C*-algebras such that there is an isomorphism between $\Gamma(A)$ and $\Gamma(B)$, in a sense made precise in Section 2, then there is an isomorphism between $A$ and $B$ which induces at the level of the invariant the given isomorphism.}\\
\begin{proof}
\indent In the previous section, Theorem 11.0.7, we have proved that any inductive limit of continuous-trace C*-algebras whose spectrum is $[0,1]$ is isomorphic to an inductive limit of special continuous-trace C*-algebras and in particular is isomorphic to an inductive limit of hereditary sub-C*-algebras with finite step function dimension function. Therefore we can use the Stevens Classification Theorem, Theorem 3.0.4, to get an isomorphism theorem for simple inductive limits of continuous-trace C*-algebras.  \\
\end{proof}
 
 \indent For completeness of the argument we reproduce the main steps of the proof of Theorem 3.0.4. For the proofs of all these steps we refer the reader to \cite{ste0}.

\subsection {Pulling back the isomorphism between inductive limits to the finite stages}
 The first step in the proof is the pulling back of the maps from the inductive limit to the finite stages. As it is shown in \cite{ste0}, at the level of the dimension range we obtain an exactly commuting diagram. At the level of the affine space, after we pull back the isomorphism between the inductive limits to the finite stages, we get an approximate intertwining between affine spaces which approximately preserves the scale.

\subsection{Existence Theorem} To construct maps between the building blocks, we shall need the Stevens Existence Theorem.
\begin{theorem}
 Let $A$ be a step hereditary sub-C*-algebra (see Definition 7.1.2 or \cite{ste0} for a definition of step-hereditary sub-C*-algebra). Let a finite subset $F$, $F$ contained in $AffT^{+}A$, and $\epsilon >0$ be given. Then there are $N \in \mathbb{N}$ and $F' \in Aff'A$ such that for $B$ a step hereditary C*-algebra and $\varphi_0 :D(A) \rightarrow D(B), \psi _{T}:AffT^{+}A \rightarrow AffT^{+}B$ maps satisfying the following conditions: \\
  1. $\varphi_0$ has multiplicity at least $N$ (recall that $D(A),D(B)$ are intervals in $\mathbb{Z};$\\
  2. dist($\varphi_{T}(F'),Aff'B) \leq 1/N;$\\
  3. $\varphi_0$ and $\varphi_T$ are approximately compatible in the sense that
  $$||\widehat{\varphi([r])}-\varphi_T(\widehat{[r]}||<1/N$$
 for any $[r] \in D(A),r \in Proj(A)$,
 there is a homomorphism $\psi : A \rightarrow B$ such that $\varphi_0 =\psi_0$ and $||(\varphi_T-\psi_T)a||<\epsilon$ for all $a \in F$.
\end{theorem}

\subsection{Uniqueness Theorem}
 At this stage we have a (not necessarily approximately commutative) diagram of algebra homomorphisms that induces an approximately commutative diagram at the level of the inductive limit systems of invariants which were obtained after pulling back the isomorphism. To make the diagram approximately commutative we need to apply the Stevens Uniqueness Theorem:

 \begin{theorem} Let A be a step hereditary sub-C*-algebra of an interval algebra and $d>0$ an integer, $F \subset A$ a finite subset and $\epsilon > 0$. Let $B$ be a step-hereditary subalgebra and $\psi, \varphi : A \rightarrow B$ with the following properties:
 1. $\varphi _{0}= \psi _{0} : D(A) \rightarrow D(B)$,\\
 2. $\psi$ and $\varphi$ have at least $\delta$-fraction eigenvalues in each $d$ consecutive subintervals of length $\frac{1}{d}$ of $[0,1]$,
 with $\hat{r_{i}}$ the functions equal to $0$ from $0$ to $\frac{i}{d}$,equal to $1$ on $[\frac{i+1}{d},1]$ and linear in between, $||(\varphi _{T}-\psi _{T})(\hat{r_{i}})||<\delta$,\\
 3. if $x_1,x_2 \in [0,1],|x_{1}-x_{2}|<\frac{3}{d}$ then $||a(x_1)-a(x_2)||<\frac{\epsilon}{2}$ for all $a\in F$. 
 Then there is an automorphism of $B$, $f$, that preserves the invariant
 and is such that 
 $$||(\psi - f\varphi)(a)||< \epsilon \;\;\;\forall a \in F$$
\end{theorem}
\begin{remark} It is a consequence of the I. Stevens's uniqueness theorem that the automorphism $f$ obtained above fixes the spectrum. It is well known, see \cite{ell4} or \cite{lan}, that an automorphisms of continuous-trace C*-algebras whose spectrum is $[0,1]$ and which fixes the spectrum is necessarily an approximately inner automorphism.  
\end{remark}
The above remark allows us to state the following Uniqueness Theorem:
 \begin{theorem} Let A be a special continuous-trace C*-algebra and $d>0$ an integer, $F \subset A$ a finite subset and $\epsilon > 0$. Let $B$ be a special continuous-trace C*-algebra and $\psi, \varphi : A \rightarrow B$ with the following properties:
 1. $\varphi _{0}= \psi _{0} : D(A) \rightarrow D(B)$,\\
 2. $\psi$ and $\varphi$ have at least $\delta$-fraction eigenvalues in each $d$ consecutive subintervals of length $\frac{1}{d}$ of $[0,1]$,
 with $\hat{r_{i}}$ the functions equal to $0$ from $0$ to $\frac{i}{d}$,equal to $1$ on $[\frac{i+1}{d},1]$ and linear in between, $||(\varphi _{T}-\psi _{T})(\hat{r_{i}})||<\delta$,\\
 3. if $x_1,x_2 \in [0,1],|x_{1}-x_{2}|<\frac{3}{d}$ then $||a(x_1)-a(x_2)||<\frac{\epsilon}{2}$ for all $a\in F$. 
 Then there is an approximately inner automorphism of $B$, $f$, such that 
 $$||(\psi - f\varphi)(a)||< \epsilon \;\;\;\forall a \in F$$
\end{theorem}
\begin{remark} Knowing that the automorphism in the uniqueness theorem is approximately inner simplifies the intertwining argument.
\end{remark}
\subsection{The Elliott approximate intertwining argument} 
\indent The last step is to observe that by applying the uniqueness theorem we can modify the maps so that the diagram approximately commutes. In other words passing to suitable subsequences of algebras, it is possible to perturb each of the homomorphisms obtained in the Existence Theorem by an inner automorphism, in such a way that the diagram becomes an approximate intertwining, in the sense of Theorem 2.1,\cite{ell3}.\\
 Therefore, by the Elliott approximate intertwining theorem (see \cite{ell3}, Theorem 2.1), the algebras $A$ and $B$ are isomorphic.
 
\subsection{Proof of Theorem 2.0.1, the isomorphism theorem for certain simple stably AI algebras}
 By Theorem 11.0.7 we get that all simple stably AI algebras which are inductive limit of continuous-trace C*-algebras with spectra the closed interval $[0,1]$ and all irreducible representations finite dimensional are inductive limits of step hereditary sub-C*-algebras. Therefore these simple stably AI algebras satisfies the hypothesis of the Theorem 3.0.3. Consequently we can apply the isomorphism Theorem 3.0.3 and prove Theorem 2.0.1. \\

\section{ The range of the invariant}

 In this section we consider the class of simple AI algebras which are not necessarily of real rank 0.
Using algebras from this class we construct a simple stably AI algebra which will exhaust the invariant.\\
\indent In order to construct a simple stably AI algebra we use the trace norm map:
\begin{definition} Let ${\mathcal H}$ be a sub-C*-algebra of a $C^*$-algebra ${\mathcal A}$. The trace norm map associated to ${\mathcal H}$ is a function $f:T^+({\mathcal H}) \rightarrow (0,\infty]$ such that $f( \tau) = || \tau|_{{\mathcal H}} ||$, $ \infty$ if $ \tau$ is unbounded.
\end{definition}
 Recall that:
\begin{definition} $T^+({\mathcal A})$ is the cone of positive trace functionals on ${\mathcal A}$ with the inherited w*-topology.
\end{definition}
\begin{remark} Let ${\mathcal A}$ be a continuous-trace C*-algebra whose spectrum is $[0,1]$. The trace norm map restricted to the extreme traces is exactly the dimension range.
\end{remark}
\begin{remark} Knowing the trace norm map is equivalent to knowing the special subset of the affine space, the scale, $Aff'()$. For instance if the trace norm map takes only the infinity value then the scale is the whole affine space $AffT^{+}()$.
\end{remark}
\begin{remark}The trace norm map must be a lower semicontinuous affine map (being a supremum of a sequence of continuous functions).
\end{remark}
\begin{remark}The dimension range can be determined using the values of the trace norm map $f$ , the simplex of tracial states $S$ and dimension group $G$. A formula for the dimension range $D$ is: 
 $$D = \{ x\in G / v(x) < f(v) ,  v \in S , v \neq 0 \}$$
\end{remark}

\begin{remark} If the algebra ${\mathcal H}$ is unital then the trace norm map does not take the infinity value.
\end{remark}

 I. Stevens constructed a hereditary sub-C*-algebra of a simple (unital) AI-algebra which is obtained as an inductive limit of hereditary subalgebras of the building blocks of the AI-algebra, and has as a trace norm map any given affine continuous function; cf.\cite{ste0}, Proposition 30.1.7. Moreover she showed that any lower semicontinuous map can be realized as a trace norm map in a special case. Our result is a generalization to the case of unbounded trace norm map when restricted to the base of the cone.\\

\noindent
{\bf Theorem 2.0.3} \textit{ Suppose that G is a simple countable dimension group, V is the  cone associated to a metrizable Choquet simplex S , $\lambda : S \rightarrow Hom^+(G,R)$ is a continuous affine map with its range dense and sending extreme rays in extreme rays, and $f:S\rightarrow (0,\infty]$ any affine lower semicontinuous map. Then $[G,(V,S), \lambda,f]$ is the Elliott invariant of some simple non-unital algebra stably AI algebra.}
\begin{proof} The proof is based on I. Stevens's proof in a special case and consists of several steps.\\
\noindent {\bf Step 0}\\
\indent We start by constructing a simple stable AI-algebra ${\mathcal A}$ with its Elliott invariant: $[(G,D),V,\lambda ]$. We know that this is possible (see \cite{ste}). By tensoring with the algebra of compact operators we may assume ${\mathcal A}$ is a simple stable AI algebra.\\
\noindent
{\bf Step 1}\\
\indent We restrict the map $f$ to the base $S$ of the cone $T^+({\mathcal A})$, where the cone $V$ is naturally identified with $T^+({\mathcal A})$. Since any lower  semicontinuous affine map $ f: S \rightarrow (0, +\infty]$ is a pointwise limit of an increasing sequence of continuous affine  positive maps, (see \cite{alf}), we can choose  $ f = \mathrm{lim} f_n$,  where $f_n$ are continuous affine and strictly positive functions. \\
\indent Moreover by considering the sequence of functions $g_n =  f_{n+1}-f_n$ if $n>1$ and $g_1 = f_1$ we get that:
 $$\sum_{n=1}^{\infty}g_n = f$$
{\bf Step 2}\\
\indent Next we use the results of Stevens (\cite{ste0}, Prop. 30.1.7), to realize each such continuous affine map $g_n$ as the norm map of a hereditary sub-C*-algebra ${\mathcal B}_n$ (which is an approximately step hereditary) of the AI-algebra ${\mathcal A}$ obtained at Step 0.\\
\indent Consider the $L^{\infty}$ direct sum $\oplus {\mathcal B}_i$ as a sub-C*-algebra of ${\mathcal A}$. The trace norm map of the sub-C*-algebra $\oplus {\mathcal B}_i$ of ${\mathcal A}$ is equal to $\sum_{i=1}^{\infty}g_n = f$.\\
\indent To see that $\oplus {\mathcal B}_i$ is a sub-C*-algebra of ${\mathcal A}$ we use that ${\mathcal A}$ is a stable C*-algebra:\\
 $$ \oplus {\mathcal B}_i = \left( \begin {array} {cccc} {\mathcal B}_1 &     & 0 \\
                                     & {\mathcal B}_2 &   \\
                                  0  &     & \ddots  \\
 \end{array} \\ \right) \subseteq {\mathcal A}\otimes \mathbb{K} \cong {\mathcal A}. $$
\indent Next we consider the hereditary sub-C*-algebra ${\mathcal H}$ generated by $\oplus {\mathcal B}_i$ inside of ${\mathcal A}$.\\
\indent In order to prove that the trace norm map of ${\mathcal H}$ is $f$ it is enough to show that the norm of a trace on $\oplus {\mathcal B}_i$ is the same as on ${\mathcal H}$.\\
\indent Let $\tau$ be a trace on any C*-algebra and $u_{\lambda}$ an approximate unit of that C*-algebra. Then it is well known, see \cite{mur}, that  $||\tau|| = lim \tau (u_{\lambda})$.\\
\indent Hence it suffices to prove that an approximate unit of the sub-C*-algebra $\oplus {\mathcal B}_i$ is still an approximate unit for the hereditary sub-C*-algebra ${\mathcal H}$.\\
\indent We shall prove first that the hereditary sub-C*-algebra generated by $\oplus {\mathcal B}_i$ coincides with the hereditary sub-C*-algebra generated by one of its  approximate units. Let $(u_{\lambda})_{\lambda}$ be an approximate unit of $\oplus {\mathcal B}_i$. Denote by ${\mathcal U}$ the hereditary sub-C*-algebra of ${\mathcal A}$ generated by $\{(u_{\lambda})_{\lambda}\}$. We want to prove that ${\mathcal U}$ is equal with ${\mathcal H}$. \\
\indent Since $(u_{\lambda})_{\lambda}$ is a subset of $\oplus {\mathcal B}_i$ we clearly have:
 $${\mathcal U} \subset {\mathcal H}.$$
\indent For the other inclusion we observe that:
 $$\forall b \in \oplus {\mathcal B}_i\; :\;b=\lim\limits_{\lambda \rightarrow \infty} u_{\lambda}bu_{\lambda}. $$
\indent Now each $u_{\lambda}bu_{\lambda}$ is an element of the hereditary sub-C*-algebra generated by $(u_{\lambda})_{\lambda}$ and hence $b \in {\mathcal U}$. Therefore $ \oplus {\mathcal B}_i \subset {\mathcal U}$ which implies that: ${\mathcal H} \subset {\mathcal U}.$\\
\indent We conclude that: ${\mathcal H} = {\mathcal U}$ and hence the trace norm map of ${\mathcal H}$ is $f$. Therefore ${\mathcal H}$ is a simple stably AI algebra with the given invariant. 
\end{proof}
\begin{remark}
 \indent The approximate unit $(u_{\lambda})_{\lambda}$ of $\oplus {\mathcal B}_i$ is still an approximate unit for the hereditary sub-C*-algebra ${\mathcal U}$. To see why this is true let us consider the sub-C*-algebra of ${\mathcal A}$ defined as follows: $ \{ h \in {\mathcal A} \;|\; h= \lim\limits_{\lambda \rightarrow \infty}u_{\lambda}h \}$.\\
\indent This sub-C*-algebra of ${\mathcal A}$ is a hereditary sub-C*-algebra. Indeed let $ 0 \leq k \leq h $ with $h= \lim\limits_{\lambda \rightarrow \infty}u_{\lambda}h$. We want to prove that $k = \lim\limits_{\lambda \rightarrow \infty}u_{\lambda}k$.\\
\indent Consider the hereditary sub-C*-algebra $\overline{h{\mathcal A}h}$ of ${\mathcal A}$ which clearly contains $h$ (because $h^2 = \lim\limits_{\lambda \rightarrow \infty }hu_{\lambda}h$). Therefore $k \in \overline{h{\mathcal A}h}$.\\
\indent Since $h= \lim\limits_{\lambda \rightarrow \infty}u_{\lambda}h$ we obtain that $u_{\lambda}$ is an approximate unit for $\overline{h{\mathcal A}h}$. In particular:
 $$k= \lim\limits_{\lambda}u_{\lambda}k$$ and hence $\{ h \in {\mathcal A} \;| h= \lim\limits_{\lambda \rightarrow \infty}u_{\lambda}h \}$ is a hereditary sub-C*-algebra of ${\mathcal A}$. Since ${\mathcal U}$ is the smallest hereditary containing $(u_{\lambda})_{\lambda}$ we get that:
 $${\mathcal U} \subset \{ h \in {\mathcal A} \;| h= \lim\limits_{\lambda \rightarrow \infty}u_{\lambda}h \} $$ 
and hence $u_{\lambda}$ is an approximate unit for ${\mathcal U}$.
\end{remark}

\section{Non-AI algebras which are stably (isomorphic) AI algebras}
In this section we present a sufficient condition on the invariant that will allow us to construct a stably AI algebra which is not an AI-algebra.\\
\indent With $[G,V, \lambda,f]$ as before we observe that for an AI-algebras with the Elliott invariant canonically isomorphic to the given invariant we have the following equality always holds:
$$f(v) = \mathrm{sup} \{v(g): g \in D \},$$ 
where D is the dimension range. This is true by simply using the fact that any AI-algebra has an approximate unit consisting of projections.\\
\indent Therefore a sufficient condition imposed on the invariant in order to get a stably isomorphic AI algebra but not an AI-algebra is :
$$f(v) \neq \mathrm{sup} \{v(g): g \in D \}. $$  
\indent This condition is also necessary. Namely assume that we have $f(v) = \mathrm{sup} \{v(g): g \in D \}$ and we have  constructed a simple ASH C*-algebra stably isomorphic to an AI algebra ${\mathcal A}$ with the invariant canonically isomorphic with the tuple $[G,V, \lambda,f]$. Consider $D = \{x \in G : v(x) < f(v), v \in S, v \neq 0\}$, where $S$ is a base of the cone $V$. For the tuple $[G,D,V,S,\lambda]$ we can build (via the range of the invariant for simple AI algebras, \cite{ste}) a simple AI-algebra ${\mathcal B}$ with the invariant naturally isomorphic with the given tuple. \\
\indent Note that the trace norm map which is defined starting from the tuple \\
$[K_{0}({\mathcal B}),D({\mathcal B}),T^{+}{\mathcal B}, \lambda_{{\mathcal B}}]$ is exactly $f$ because of the equality: 
$$f(v) = \mathrm{sup} \{v(g): g \in D \}$$
  and ${\mathcal B}$ is an AI-algebra.\\
\indent It is clear that ${\mathcal B}$ is a stably AI algebras and hence by the isomorphism Theorem 2.0.1 for the class of simple ASH C*-algebras stably isomorphic AI algebras we conclude that ${\mathcal A}$ isomorphic to ${\mathcal B}$. Hence ${\mathcal A}$ is a simple AI algebras as desired and we have proved the following theorem:  
\begin{theorem} Let ${\mathcal A}$ be a simple C*-algebra stably AI algebra which is an inductive limit of continuous-trace C*-algebras whose spectrum is homeomorphic to $[0,1]$ and all irreducible representations are finite dimensional. A necessary and sufficient condition for ${\mathcal A}$ to be a simple AI algebra is:
$$f(v) = \mathrm{sup} \{v(g): g \in D \}. $$  
\end{theorem}

\section{ The size of the class of simple stably AI algebras is much larger then the size of the class of simple AI algebras}  To see this we consider the simple AI algebra necessarily not of real rank zero with scaled dimension group $(\mathbb{Q},\mathbb{Q}_{+})$ and cone of positive trace functionals a 2-dimensional cone; see \cite {ste}. Then the set of possible stably AI algebras, or equivalently the set of possible trace norm maps, may be represented as the extended affine space shown in the following schematic diagram:

\begin{center}
\includegraphics [height =3cm]{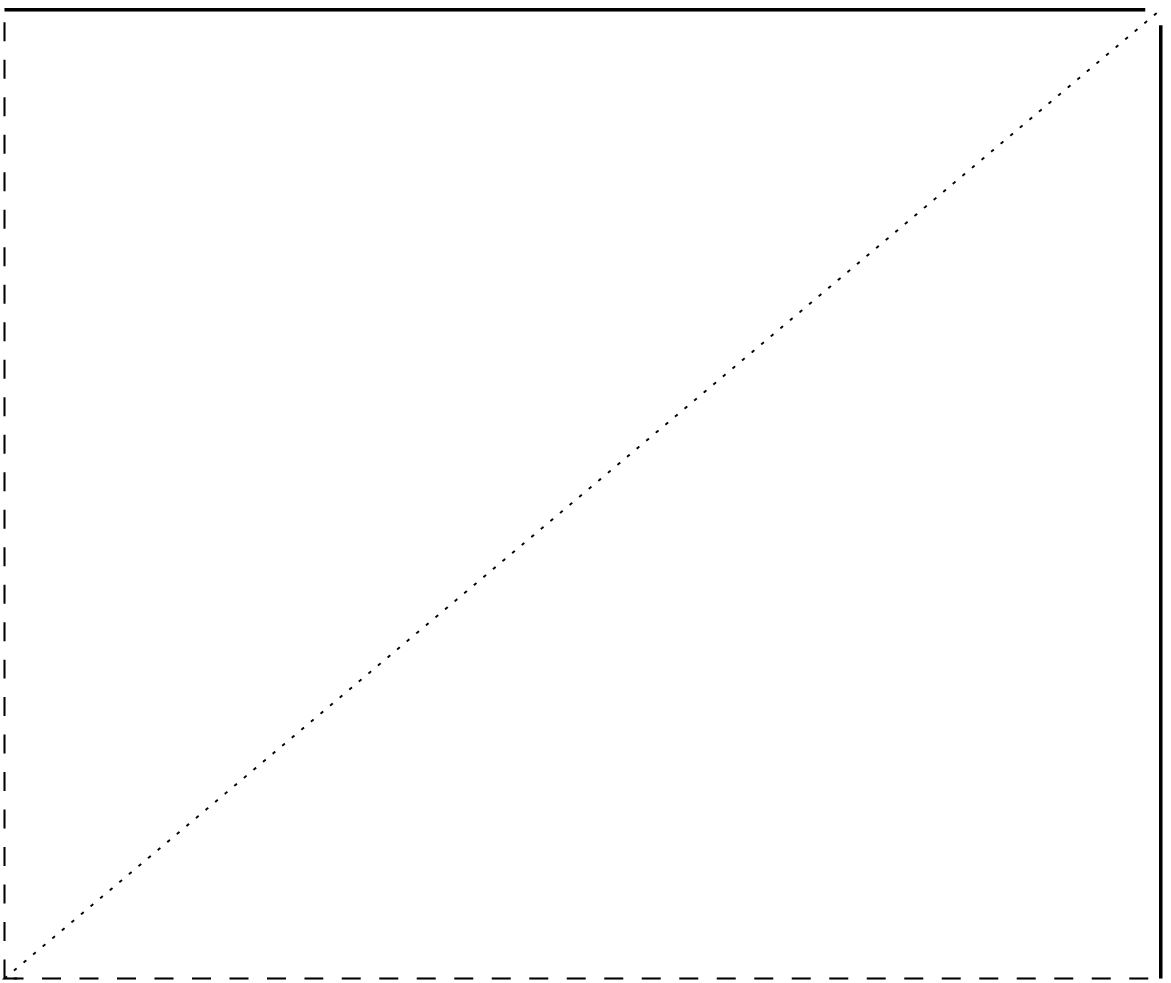}
\end{center}

\indent Each off-diagonal point in the diagram is the trace norm map of one of I. Stevens's algebras. The boundary points of the first quadrant are removed (dotted lines) and the points with infinite coordinates are allowed. The dimension range is embedded in a canonical way in the extended affine space as the main diagonal consisting of the points with rational coordinates.\\
\indent The two bold lines represent the cases of stably AI algebras with unbounded trace norm map (points on these two lines have at least one coordinate infinity).\\
\indent If the point is off the diagonal and in the first quadrant, by Theorem 13.0.9 we get that the corresponding stably AI algebra is an ASH algebra which is not AH. It is clear that the size of the set of points off the diagonal is much larger then the size of the set of points on the diagonal.(For instance in terms of the Lebesgue measure.)\\
\indent This picture shows that the class of simple AI algebras sits inside the class of stably AI algebras as the main diagonal sits inside the first quadrant. (Note especially the small dots which correspond to both unital and non-unital AI algebras: recall the case of stably UHF algebras!)\\
\indent The algebras in the preceding example are not of real rank zero.
\bibliographystyle{amsplain}

\end{document}